\documentclass{ws-aa}
\usepackage{amssymb,wasysym}
\newcommand\R{{\mathbb R}}

\newcommand\Dis{\displaystyle}
\newcommand\fa{\frac}
\newcommand\veps{\varepsilon}

\begin{document}

%\markboth{Authors' Names}{Instructions for Typing Manuscripts (Paper's Title)}

%%%%%%%%%%%%%%%%%%%%% Publisher's Area please ignore %%%%%%%%%%%%%%%
%
%\catchline{}{}{}{}{}
%
%%%%%%%%%%%%%%%%%%%%%%%%%%%%%%%%%%%%%%%%%%%%%%%%%%%%%%%%%%%%%%%%%%%%
\title{Semi classical measures and Maxwell's system\\}

\author{Hassan TAHA}
\address{Universit\'e d'Orl\'eans \\
MAPMO-UMR 6628, BP 6759\\
45067 ORLEANS CEDEX-France\\ 
hassan.taha@labomath.univ-orleans.fr}
\maketitle

\begin{abstract}
\noindent We are interested in the homogenization of 
energy like quantities for electromagnetic waves in the high frequency limit
for Maxwell's equations with various boundary conditions. 
We use a scaled variant of H-measures known as semi classical measures or 
Wigner measures.\par
\noindent Firstly, we consider this system in the half space 
of $\R^3$ in the time harmonic and with conductor boundary condition at 
the flat boundary $x_3=0$. Secondly we consider the same 
system but with Calderon boundary condition. Thirdly, we consider 
this system in the curved interface case.\par   
\end{abstract}


\noindent \section{\bf Introduction}\noindent In this work, we are interested in the homogenization of energy quantities for electromagnetic waves in the high frequency limit, and more particularly for Maxwell's equations. Our interest is also in dealing with interactions with various boundary conditions. For this purpose, we use a scaled variant of H-measures (see L.Tartar or P.G\'erard), known as semi classical measures or Wigner measures, introduced in \cite{Gér2:Mar}, \cite{Geoleon}, \cite{Geoleonkel}.\par\noindent One of the most important predictions of Maxwell's equations is the existence of electromagnetic waves which can transport energy.\par \noindent For this purpose, the Theory of Radiative Transfer was originally developed to describe how light energy propagates throught a turbulent atmosphere. This theory can applied to representative problems involving reflection, transmission, and diffraction in both homogeneous and inhomogeneous media.\par
\noindent Justification of this theory in high frequency limit, as well as for other waves equations, can be 
found for a deterministic medium in the works of P.Gerard \cite{Gér1}, \cite{Gér3}, and C.Kammerer \cite{caclotid} 
as well as by P.L.Lions and T.Paul \cite{Lipa} and L.Miller \cite{lmul} and 
G.Papanicolaou \cite{GeoleonkelGui}.\par \noindent Our purpose in this paper is to describe this energy propagation for Maxwell's system, coupled with 
various boundary conditions, and with a typical scale which is played here by the frequency.\par
\smallskip
\noindent We shall consider Maxwell's system, with electric permeability $ \ddot{\epsilon} $, conductivity 
$\ddot{\sigma} $ and magnetic susceptibility $\ddot{\eta} $, in the half space $(x^{3}\geq 0)$ of $\R^{3}$, with the courant variable $x =x^1,x^2,x^3$. These quantities are $3\times 3$ matrix valued functions of $x$. This system is given by the following equations\par
\smallskip

\smallskip
\begin{equation}\label{maxwel}\left\{\begin{array}{ccccccc}i) &\partial_{t}D^{\Dis\varepsilon} (x,t) +J^{\Dis\varepsilon} (x,t) &= &\mbox{rot} H^{\Dis\varepsilon} (x,t) +F^{\Dis\varepsilon}&,\\[0.3cm]     ii) &\partial_{t} B^{\Dis\varepsilon} (x,t) &= &-\mbox{rot} E^{\Dis\varepsilon} (x,t) +G^{\Dis\varepsilon} (x,t)&,\\[0.3cm]     iii) &\mbox{div} B^{\Dis\varepsilon} (x,t)& =&0&,\\[0.3cm] iv) &\mbox{div} D^{\Dis\varepsilon} (x,t)& =&\varrho^{\Dis\varepsilon} (x,t)&   \end{array}\right.\end{equation}
\noindent where $t\in(0,T)$, and $ E^{\;\Dis\veps}\;,H^{\;\Dis\veps}\;,D^{\;\Dis\veps}\;,J^{\;\Dis\veps}$ 
and $B^{\;\Dis\veps}$ are the electric, magnetic, induced electric, current density and induced magnetic fields, 
respectively. Morever, $ \rho^{\;\Dis\veps}$ is the charge density (a function uniformly bounded in $L^2 (\R^3 )$, and where $F^{\;\Dis\veps}$, $ G^{\;\Dis\veps}\in L^{2}(\R^{3})^{3}$ are given.\par\smallskip
\smallskip
\noindent We complete this system by the following constitutive relations
\begin{equation}\label{diss}\left\{\begin{array}{ccccccc}    1) &\Dis D^{\Dis\veps}(x,t)&=&\ddot{\epsilon}(x) E^{\Dis\veps}(x,t)& ,\\[0.3cm] 2) &\Dis J^{\Dis\veps}(x,t)&=&\ddot{\sigma}(x)E^{\Dis\veps}(x,t)&,\\[0.3cm] 3) &\Dis B^{\Dis\veps}(x,t)&=&\ddot{\eta}(x)H^{\Dis\veps}(x,t)&\;.\\[0.3cm]      \end{array}\right.\end{equation}

\noindent We shall only be interested in time harmonic solutions of this system and in the high frequency limit. For that purpose, we look for solutions in the form
\begin{equation}\label{maxx}\left\{\begin{array}{ccccccc}
&\Dis D^{\Dis\veps}(x,t)&=& D^{\Dis\veps}(x) \Re\{\exp \fa{i\omega t}{\Dis\veps}\}&,\\[0.3cm] 
&H^{\Dis\veps}(x,t)&=&H^{\Dis\veps}(x) \Re\{\exp \fa{i\omega t}{\Dis\veps}\}&,\\[0.3cm] 
&\Dis J^{\Dis\veps}(x,t)&=& J^{\Dis\veps}(x) \Re\{\exp \fa{i\omega t}{\Dis\veps}\}&,\\[0.3cm] 
&\Dis B^{\Dis\veps}(x,t)&=& B^{\Dis\veps}(x) \Re\{\exp \fa{i\omega t}{\Dis\veps}\}&,\\[0.3cm] 
&E^{ \Dis\veps}(x,t)&=&E^{\Dis\veps}(x) \Re\{\exp \fa{i\omega t}{\Dis\veps}\}&,\\[0.3cm] 
\end{array}\right.\end{equation}
\noindent where $\omega$ is the given fixed frequence, that we assume different from $0$. Note that we use the same letters on both sides of the above equations to simplify notations.\par 
\smallskip\noindent In this work, we assume that the matrix $\ddot{\epsilon}$ $\;,\ddot{\eta}$$\;,\ddot{\sigma} $, are $3\times 3$ "scalar" matrix valued functions given by
\begin{equation}\label{3.25}
\ddot{\epsilon}= \epsilon  ({\bf Id})_{3\times 3} \equiv \left(               \begin{array}{cccccccccccccccccccccc}                  \epsilon(x)&   0  &   0\\[0.3cm]                     0&             \epsilon(x)  &   0\\[0.3cm]                     0&             0  &   \epsilon(x)\\[0.3cm]                     \end{array}                 \right)\end{equation} 
\noindent and

\begin{equation}\label{3.26}\ddot{\eta}=\eta ({\bf Id})_{3\times 3} \equiv  \left(               \begin{array}{cccccccccccccccccccccc}                  \eta (x)&    0  &    0\\[0.3cm]                        0 &           \eta (x)   &   0\\[0.3cm]                        0 &           0   &  \eta (x) \\[0.3cm]                        \end{array}                 \right)\;,\ddot{\sigma}=\sigma ({\bf Id})_{3\times 3} \equiv \left(               \begin{array}{cccccccccccccccccccccc}                  \sigma(x)&   0   &   0\\[0.3cm]                       0&              \sigma(x)   &   0\\[0.3cm]                       0&              0   &   \sigma(x)\\[0.3cm]                        \end{array}                 \right)\end{equation}

\noindent where $ \epsilon\;,\eta\;,\sigma $, are smooth (scalar) functions in $C^{1}(\R^3) $. This usual assumption could be certainly relaxed, but at the expense of much more complex spectral calculus.\par 
\noindent With the above notations, the time harmonic form of Maxwell equations are then
\begin{equation}\label{rajout1}\left\{\begin{array}{ccccccc} &\mbox{rot} E^{\Dis\varepsilon}-i\omega {\bf \eta} H^{\Dis\varepsilon}=F^{\Dis\varepsilon}&,\\[0.3cm] &\mbox{rot} H^{\Dis\varepsilon}+i\omega {\bf \epsilon} E^{\Dis\varepsilon}=G^{\Dis\varepsilon}.&\end{array}\right.\end{equation}
\noindent Note that we have not written the third and fourth equations appearing in system (\ref{maxwel}), since in fact we assume that the right hand sides of (\ref{rajout1}) do satisfy the usual compatibility conditions.\par
\noindent Set$$u^{\Dis\veps}= \left(                \begin{array}{cccccccccccccccccccccc}                        E^{\Dis\veps}\\                     H^{\Dis\veps}\\                 \end{array}\right )=                  \left(                 \begin{array}{cccccccccccccccccccccc}                        E_1^{\Dis\veps}\\                     E_2^{\Dis\veps}\\                    E_3^{\Dis\veps}\\                     H_1^{\Dis\veps}\\                     H_2^{\Dis\veps}\\                     H_3^{\Dis\veps}\\   \end{array}    \right) \;,$$
\vspace{0.4cm}
\begin{equation}\label{3.rr5} A^{\;0} =\left(                  \begin{array}{cccccccccccccccccccccc}                      \bf{\epsilon} \bf{Id}&  \bf{0}  \\                     \bf{0}      &  \bf{\eta} \bf{Id} \\
  \end{array}    \right)\end{equation}  
\noindent and

\begin{equation}\label{elke}  A^{\;1} =\left(                  \begin{array}{cccccccccccccccccccccc}                        \bf{0} &  \bf{Q_1}^{t}  \\                     \bf{Q_1}    &  \bf{0}       \\  \end{array}    \right)\;, A^{\;2} =\left(                  \begin{array}{cccccccccccccccccccccc}                        \bf{0} &  \bf{Q_2}^{t}  \\                     \bf{Q_2}    &  \bf{0}     \\  \end{array}    \right)\;, A^{\;3} =\left(                 \begin{array}{cccccccccccccccccccccc}                        \bf{0} &  \bf{Q_3}^{t}  \\                     \bf{Q_3}    &  \bf{0}   \\  \end{array}    \right)\;\end{equation}    
\vspace{0.3cm}
\noindent where the constant antisymmetric matrices ${\bf Q_k} \;, 1\leq k\leq 3$ are given by
\begin{equation} \bf{Q_1} =\left(                  \begin{array}{cccccccccccccccccccccc}                        0 &  0  & 0 \\                     0 &  0  & -1  \\                     0 &  1  & 0   \\   \end{array}    \right)\;,\bf{Q_2} =\left(                  \begin{array}{cccccccccccccccccccccc}                        0 &  0  & 1 \\                     0 &  0  & 0  \\                    -1 &  0  & 0   \\   \end{array}    \right)\;,\bf{Q_3} =\left(                  \begin{array}{cccccccccccccccccccccc}                        0 &  -1  &  0 \\                     1 &   0  &  0  \\                     0 &   0  &  0   \\   \end{array}  
  \right)\;.\end{equation}
\vspace{0.3cm}
\noindent Above, the matrix $ C $ is given by
\begin{equation}\label{3.522}C= \left(   \begin{array}{cccccccccccccccccccccc}       \bf{\sigma}\bf{Id}& \bf{0}\\     \bf{0}     & \bf{0} \\   \end{array}    \right)\end{equation}
\vspace{0.3cm}
\noindent while the right hand side is $ f^{\Dis\veps}= \left(                \begin{array}{cccccccccccccccccccccc}                        {F}^{\Dis\veps}\\                     {G}^{\Dis\veps}\\                 \end{array}\right )$.\par

\noindent Assumed uniform boundedness and symmetry of the permeability and susceptibility tensors show that 
system (\ref{rajout1}) is a symmetric hyperbolic system as follows   
\begin{equation}\label{3.4}\fa{i\omega}{\displaystyle\veps} A^{0}+\Dis\sum_{j=1}^{3} A^{j} \frac{\partial u^{\displaystyle\veps}}{\partial x_{j}}+C u^{\displaystyle\veps}= f^{\displaystyle\veps}.\end{equation}    
\smallskip
\noindent As a first boundary value problem, we shall consider system (\ref{rajout1}) or equivalently system (\ref{3.4}), posed in a domain, that we choose to be the upper half plane $\R^3_+ = \lbrace x, x_3 \geq 0 \rbrace$, with a perfect conductor boundary condition at the flat boundary $x_3 =0$, i.e\par
\smallskip
\smallskip
\begin{equation}\hspace{2cm} \vec{n}^+\wedge E^{\Dis\veps}=0 \; \hspace{1.5cm} \mbox{on} \hspace{1cm} x_{3}=0 \end{equation}\noindent where $\vec{n}^+=\left(                \begin{array}{cccccccccccccccccccccc}                        0\\                     0\\
                     -1\\                 \end{array}\right)$ is the unit outward normal vector to $\R^3_+$. Note that this domain is not bounded, but this is unimportant since we will localize all our functions.\par 
\smallskip
\noindent The second problem dealt with in this paper will be a transmission problem. To simplify the exposition, we will consider a medium, made of two parts: $\R^3_+\equiv\lbrace x\in \R^{3}, x_3\geq 0\rbrace$ will be the exterior medium, while $\R^3_-\equiv \lbrace x\in\R^{3}, x_3\leq 0\rbrace$ will be the interior one, each caracterized by distinct electromagnetic coefficients.\par
\noindent We suppose that our electromagnetic field is created by an incident wave 
$u^{inc}= \left(                \begin{array}{cccccccccccccccccccccc}                        E^{inc}\\                     H^{inc}\\                 \end{array}\right)$.\par
\noindent In $\R^3_+$, we consider the following exterior problem, characterized by the dielectric 
coefficients $({\bf \epsilon}^{ext}(x)\;, {\bf \eta}^{ext}(x))$ belonging to $C^{1}(\R^{3})$, 
and scalar valued, see \cite{ngde} \par 

\begin{equation}\label{exte}\left\{\begin{array}{ccccccc} &\mbox{rot} E^{ext,\Dis\varepsilon}-i\omega {\bf \eta}^{ext} H^{ext,\Dis\varepsilon}=0&,\\[0.3cm] &\mbox{rot} H^{ext,\Dis\varepsilon}+i\omega {\bf \epsilon}^{ext} E^{ext,\Dis\varepsilon}=0&,\\[0.3cm] &\left|\sqrt{{\bf \epsilon}^{ext}}E^{ext,\Dis\varepsilon}-\sqrt{{\bf \eta}^{ext}}H^{ext,\Dis\varepsilon}\wedge n^+\right|\leq\Dis\fa{c}{r^{2}}.&\end{array}\right.\end{equation}
\noindent Here $E^{ext,\Dis\veps}\;,H^{ext,\Dis\varepsilon}$ are the so called exterior fields, $r=|x|$, $\vec{x}=(x_1,x_2,x_3)$ and the third equation is 
the classical Silver Muller radiation condition, see for more details \cite{cess}, \cite{ngde}, with $n^+$ being the unit outward normal vector to $\R^3_+$.\par   
\smallskip
\smallskip

\noindent In $\R^{3}_-$, we consider the following interior problem, which is characterized by the dielectric coefficients $({\bf \epsilon}^{int}(x)\;, {\bf \eta}^{int}(x))$ belonging to $C^{1}(\R^{3})$, and scalar valued, see \cite{ngde}  
\begin{equation}\label{inte}\left\{\begin{array}{ccccccc} &\mbox{rot} E^{int, \Dis\varepsilon}-i\omega {\bf \eta}^{int} H^{int,\Dis\varepsilon}=0&,\\[0.3cm] &\mbox{rot} H^{int,\Dis\varepsilon}+i\omega {\bf \epsilon}^{int} E^{int,\Dis\varepsilon}=0.&\end{array}\right.\end{equation}
\noindent We impose the following boundary conditions (Calderon condition)
\begin{equation}\label{bound}\left\{\begin{array}{ccccccc} &E_{\Dis\varepsilon}^{int}\wedge n^--(E_{\Dis\varepsilon}^{ext}+ E^{inc})\wedge n^-=0& &\mbox{on} &x_{3}=0\;\;,\\[0.3cm] &H_{\Dis\varepsilon}^{int}\wedge n^--(H_{\Dis\varepsilon}^{ext}+H^{inc})\wedge n^-=0& &\mbox{on} &x_{3}=0\end{array}\right.\end{equation}where $n^-$ is the unit outward normal vector to $\R^3_-$.\par
\noindent $E^{int,\Dis\veps}\;,H^{int,\Dis\varepsilon}$ are the so called interior fields. Note that there is no condition at infinity in the interior problem, mainly because we have assumed intuitively a localization near $x_3 =0$.\par
\smallskip
\smallskip

\noindent In the third and final part, we generalise these two cases, and we study the curved interface case, where 
the plane $x_3 =0$ is now replaced by a curved interface, in the spirit of the work of G\'erard and Leichtman
\cite{Géreric}.\par 
%\noindent {\bf Ici, rajouter le cadre et un petit resume, etant donne %qu'il n'y a rien de mis; sinon on ne comprend pas le theoreme!!}\par
\noindent More precisely, we consider Maxwell's system (\ref{rajout1}) given above the surface given by $\Gamma:\ x_{3}=\phi(x^{'})$, where $x^{'}=(x_{1},x_{2})$, and $\phi\in W^{2}(\R^{2},\R)$ is a scalar function.\par\noindent We consider this system in time harmonic form, in the high frequency limit, and we consider a perfect boundary condition on $\Gamma $.\par
\smallskip\noindent  For each of the above cases, we shall study propagation of energy like quantities, using the framework of semi classical measures. Basic facts about these tools are recalled in Section II, refering the reader for more details to \cite{Gér2}, \cite{Geoleon}.\par\noindent  Then in Section III, we consider the above cases of Maxwell's equations, with differents 
boundary conditions, and in particuliar, we prove therein the following results
\noindent\begin{theorem}{\bf Perfect conductor case} \label{hsdaa}\noindent Consider time harmonic Maxwell's system in the half space $x_3\geq 0$ with a perfect boundary condition, written in the form (\ref{3.4}), with solution vector $u^{\Dis\veps}$. 
Let $\theta(x)$ be a test function with compact support that is equal to one on a compact set $K\subset \R^{3}$.
Let $ u^{\theta ,\Dis\veps} = \theta u^{\Dis\veps}$ be uniformly bounded in $ L^{2}(\R^{3})$, with (up to a subsequnce) an associated semiclassical measure $ \ddot{\mu}$.
\smallskip
\noindent Then the semi classical  measure 
$ \ddot{\mu}$ is supported on the set ($x\in \ Supp\ \theta , k \in \R^3$)\begin{equation}\label{support}U= \left\{(x,k)\;,\;\omega_{+}=\omega \right\}\cup \left\{(x,k)\;,\; \omega_{-}=\omega \right\}\end{equation}
\smallskip\noindent where $v(x)=\Dis\fa{1}{\sqrt{\epsilon(x)\eta(x)}}\;$ is the propagation speed. Above 
$\omega_{0}= \omega_0 (x,k) =0$, $\omega_{+}= \omega_+ (x,k)=v(x)|k|$, $\omega_{-}=\omega_- (x,k)=-v(x)|k|$ 
are the eigenvalues (of constant multiplicity two) of the dispersion matrix $L(x,k)=\Dis\sum_{j=1}^{3}
(A^{0})^{-1} k_{j} A^{j}$.\par
\smallskip
\noindent The semi classical measure $\ddot{\mu}(x,k)$ has the form 
\begin{equation}\label{theo-decomposition}
\left\{
\begin{array}{ccccccc}\ddot{\mu}(x,k)= \mu_{+}^{1}(x,k) b_{+}^{1}(x,k)\otimes b_{+}^{1*}(x,k)+
\mu_{+}^{2}(x,k) b_{+}^{2}(x,k)\otimes b_{+}^{2*}(x,k)\vspace{0.4cm}\\
+\mu_{-}^{1}(x,k) b_{-}^{1}(x,k)\otimes b_{-}^{1*}(x,k)+
\mu_{-}^{2}(x,k) b_{-}^{2}(x,k)\otimes b_{-}^{2*}(x,k)
\end{array}  
\right.
\end{equation}
\smallskip\noindent where $\mu_{+}^{1}\;,\mu_{+}^{2}$ are two scalar positive measures supported on the set $\left\{(x,k)\;,\;\omega_{+}=\omega \right\}$ and $\mu_{-}^{1}\;, \mu_{-}^{2}$, are two scalar positive measures supported on the set $\left\{(x,k)\;,\; \omega_{-}=\omega \right\}$. $b_{+}^{1}\;, 
b_{+}^{2}\;$ (resp. $ b_{-}^{1}\;,  b_{-}^{2}$) are two (normalized) eigenvectors of the matrix $L(x,k)$, 
corresponding to the eigenvalue $\omega_{+}\;$ (resp. $\omega_{-}$).\par
\noindent Furthermore, the scalar measure $\mu_{+}^{1}$ satisfies the following transport equation
\smallskip
\begin{equation}\label{theo-transport}\nabla_{k}\omega_{+}.\nabla_{x}\mu_{+}^{1}- \nabla_{x} \omega_{+}.\nabla_{k}\mu_{+}^{1} = vk_3\hat k_3[\nu^1_{\alpha +}T_1\delta_{k_3= k^-_3}  
+ \nu^1_{\beta +} T_1\delta_{k_3= k^+_3} ]\delta_{x_3 =0}\end{equation}
\smallskip\noindent where $\nu_{\alpha +}^{1}\;, \nu_{\beta +}^{1}$ are scalar positive measures associated with the semiclassical measure $\ddot\nu$
corresponding to the boundary term $u^{\Dis\varepsilon ,\theta}(x',0)$, with $x'=(x_1,x_2)$ and $\hat k =k/\mid k\mid$. The 
wave vector $k^{\pm}(k')=(k', k_{3}^{\pm})$ is defined by 
$$ k_{3}^{\pm}(x',0)=\pm\sqrt{\fa{\omega^{2}}{v(x',0)^{2}}-k'^{2}} , \ k' =(k_1, k_2)\;.$$ 
\noindent Finally, we have denoted by $T_1$ the operator defined as follows: for all smooth function $a(x,k)$ let the unique decomposition of $a$ given by
$$a(x,k) =a_0 (x,k') + a_1 (x,k')k_3 + a_2(x,k) (v\mid k\mid -\omega ).$$
Then we set $T_i (a) =a_i\;, i=0,1,2$.\par
\noindent Similar results hold true for the other scalar semi-classical measures.
\end{theorem}
\smallskip
\noindent\begin{theorem} \label{helen2}\noindent {\bf Calderon boundary condition case} \noindent Using the same framework as in Theorem \ref{hsdaa}, but with Calderon boundary condition, the associated semi classical measure $\ddot{\mu}^{ext}(x,k)$, corresponding to the exterior part,
is supported on the set
\begin{equation}\label{theo-support-calderon}U= \left\{(x,k)\;,\;\omega_{+}^{ext}=\omega \right\}\cup \left\{(x,k)\;,\; \omega_{-}^{ext}=\omega \right\}\;.\end{equation}
\smallskip\noindent Furthermore, it has the form
\begin{equation}\label{theo-decomposition-exterieur}
\left\{
\begin{array}{ccccccc}\ddot{\mu}^{ext}(x,k)= \mu_{+}^{ext,1}(x,k) b_{+}^{ext,1}(x,k)\otimes b_{+}^{ext,1*}(x,k)+
\mu_{+}^{ext,2}(x,k) b_{+}^{ext,2}(x,k)\vspace{0.4cm}\\
\otimes b_{+}^{ext,2*}(x,k)+\mu_{-}^{ext,1}(x,k) b_{-}^{ext,1}(x,k)\otimes b_{-}^{ext,1*}(x,k)+
\mu_{-}^{ext,2}(x,k) b_{-}^{ext,2}(x,k)\vspace{0.4cm}\\
\otimes b_{-}^{ext,2*}(x,k)
\end{array}  
\right.
\end{equation}
\smallskip
\noindent where $\mu_{+}^{ext,1}$, $\mu_{+}^{ext,2}$ are two scalar positive measures supported on the set $\left\{(x,k)\;,\;\omega_{+}^{ext}=\omega \right\}$, and $\mu_{-}^{ext,1}$, $\mu_{-}^{ext,2}$, are two scalar positive measures 
supported on the set $\left\{(x,k)\;,\; \omega_{-}^{ext}=\omega \right\}$. $b_{+}^{ext, 1}\;, 
b_{+}^{ext,2}\;$ (resp. $ b_{-}^{ext,1}\;,  b_{-}^{ext,2}$) are two eigenvectors of the exterior dispersion matrix 
$$L^{ext}(x,k)=\Dis\sum_{j=1}^{3}(A^{ext,0})^{-1} k_{j} A^{j},$$ corresponding to the eigenvalue 
$\omega_{+}^{ext}\;$ (resp. $ \omega_{-}^{ext}$).\par
\smallskip
\noindent The scalar transport equation, for the first scalar positive measure $\mu^{ext,1}_+$ is given by 
\begin{equation}\label{theo-transport-exterieur}\nabla_{k}\omega_{+}^{ext}.\nabla_{x}\mu^{ext,1}_+- \nabla_{x} \omega_{+}^{ext}.\nabla_{k}\mu^{ext,1}_+ =\\[0.3cm]v^{ext}\hat{k}_{3}[\nu_{+,\alpha}^{ext,1}\delta_{k_{3}=k_{3}^{ext,-}}+\nu_{+,\beta}^{ext,1}\delta_{k_{3}=k_{3}^{ext,+}}]\delta_{x_3=0}\;.\end{equation}
\smallskip\noindent Above $\nu_{+,\alpha}^{ext,1}\;, \nu_{+,\beta}^{ext,1}$ are scalar positive measures associated with the semiclassical measure $\ddot\nu^{ext}$
corresponding to the boundary term $u^{ext,\Dis\varepsilon ,\theta}(x',0)$. The wave vector $k^{ext,\pm}(k')=(k', k_{3}^{ext,\pm})$ is defined by 
$$ k_{3}^{ext,\pm}(x',0)=\pm\sqrt{\fa{\omega^{2}}{v^{ext}(x',0)^{2}}-k'^{2}}$$ 
\noindent and $v ^{ext}(x)=\Dis\fa{1}{\sqrt{\epsilon^{ext}(x)\eta^{ext}(x)}}$ is the propagation speed for the exterior problem.\par\noindent Similar results hold true for the other scalar positive measures.\par
\smallskip
\noindent Similarly, for the interior problem ( $x_{3}\leq 0$), with the following interior dispersion matrix$$L^{int}(x,k)=\Dis\sum_{j=1}^{3}(A^{int,0})^{-1}k_{j}A^{j}$$ 
\noindent the corresponding semi classical measure $\ddot{\mu}^{int}(x,k)$, is supported on the set
\begin{equation}\label{theo-support-interieur}U= \left\{(x,k)\;,\;\omega_{+}^{int}=\omega \right\}\cup \left\{(x,k)\;,\; \omega_{-}^{int}=\omega \right\}\;.\end{equation}
\noindent Furthermore, it has the form
\begin{equation}\label{theo-decomposition-interieur}
\left\{
\begin{array}{ccccccc}\ddot{\mu}^{int}(x,k)= \mu_{+}^{int,1}(x,k) b_{+}^{int,1}(x,k)\otimes b_{+}^{int,1*}(x,k)+
\mu_{+}^{int,2}(x,k) b_{+}^{int,2}(x,k)\vspace{0.2cm}\\
\otimes b_{+}^{int,2*}(x,k) +\mu_{-}^{int,1}(x,k) b_{-}^{int,1}(x,k)\otimes b_{-}^{int,1*}(x,k)+
\mu_{-}^{int,2}(x,k) b_{-}^{int,2}(x,k)\vspace{0.2cm}\\
\otimes b_{-}^{int,2*}(x,k)
\end{array}  
\right.
\end{equation}
\smallskip
\noindent where $\mu_{+}^{int,1}\;,\mu_{+}^{int,2}$ are two scalar positive measures 
supported on the set $\left\{(x,k)\;,\;\omega_{+}^{int}=\omega \right\}$ and 
$\mu_{-}^{int,1}\;, \mu_{-}^{int,2}$, are two scalar positive measures 
supported on the set $\left\{(x,k)\;,\; \omega_{-}^{int}=\omega \right\}$. $b_{+}^{int,1}\;, 
b_{+}^{int,2}\;$ (resp. $ b_{-}^{int,1}\;,  b_{-}^{int,2}$) are two eigenvectors of the matrix $L^{int}(x,k)$ given above, corresponding to the eigenvalue $\omega_{+}^{int}\;$ (resp.$ \omega_{-}^{int}$).\par
\noindent The scalar transport equation for the first positive measure $\mu^{int,1}_+$ is given by
\begin{equation}\label{theo-transport-interieur}\nabla_{k}\omega_{+}^{int}.\nabla_{x}\mu^{int,1}_+- \nabla_{x} \omega_{+}^{int}.\nabla_{k}\mu^{int,1}_+ =\\[0.3cm]v^{ext}\hat{k}_{3}[\nu_{+,\alpha}^{int,1}\delta_{k_{3}=k_{3}^{int,-}}+\nu_{+,\beta}^{int,1}\delta_{k_{3}=k_{3}^{int,+}}]\delta_{x_3=0}\end{equation}%\end{rem} 
\smallskip
\noindent where $\nu_{+,\alpha}^{int,1}\;, \nu_{+,\beta}^{int,1}$ are scalar measures associated to the semiclassical measure $\ddot \nu^{int}$
corresponding to the boundary term $u^{int,\Dis\varepsilon ,\theta}(x',0)$, and the 
tangential vector $k^{int,\pm}(k')=(k', k_{3}^{int,\pm})$ is defined by 
$$ k_{3}^{int,\pm}(x',0)=\pm\sqrt{\fa{\omega^{2}}{v^{int}(x',0)^{2}}-k'^{2}}$$ 
\noindent where $v^{int}(x)=\Dis\fa{1}{\sqrt{\epsilon^{int}(x)\eta^{int}(x)}}$ is 
the propagation speed for the interior problem.\par\noindent Similar results hold true for the other scalar positive measures. Finally, setting
$$M=\left(       \begin{array}{cccccccccccccccccccccc}  
      1 &  0  & 0 & 0 &0 &0 \\       0 &  1  & 0 & 0 &0 &0 \\       0 &  0  & 0 & 0 &0 &0 \\
       0 &  0  & 0 & 1 &0 &0 \\
       0 &  0  & 0 & 0 &1 &0 \\
       0 &  0  & 0 & 0 &0 &0 \\
    \end{array}      \right)$$
we have the following relation
$$\ddot \nu^{int} = M\ddot \nu^{ext}.$$

\end{theorem}
\smallskip

\smallskip
\noindent By adapting the proofs of the above two main theorems, we are also able to deal with the curved interface case. We sketch the proof at the end of Section III.\par
\noindent \section{Prerequesites on semi classical measures}\par\noindent In this section, we recall some properties of semi classical measures which are useful 
in the analysis of high frequency propagation problems. For more details, we refer to \cite{Gér2}, \cite{caclotid}, 
\cite{Lipa}, \cite{Geoleonkel}, \cite{Geoleon}.\par\smallskip
%\begin{Def} {\bf Wigner measure and somes properties} \par%\end{Def}
%\begin{Def}
\noindent Let $f:\R^{d}\longrightarrow \R^{n}$ be in $L^2(\R^{d})^n$. Its (unscaled) semi classical transform is then defined as
\begin{equation}W(x,k)= \Dis\fa{1}{(2\pi)^{d}} \Dis\int_{\R^{d}} e^{iy.k} f(x-y/2) \otimes f^{*}(x+y/2) dy\;.\end{equation}%\end{Def} 
\noindent Its scalar semi classical transform is $w(x,k)=Tr(W(x,k))$. The function $f$ can be scalar ($n=1$), 
or vector-valued ($f^{*}$ denotes the transposed conjugated of the vector $f$). In the latter case its semi 
classical transform is an hermitian $ n\times n $ matrix.\par
\noindent We want to consider the semi classical transform of high frequency waves, i.e of functions 
$f^{\Dis\veps}(x)$ which are oscillating on a given scale $\Dis\veps $, such that $\Dis\veps\to 0$. Our 
exposition follows the ideas of P. Gerard \cite{tar1}, \cite{Geoleon}. Therefore, we consider the rescaled 
semi classical transform, at the scale $\varepsilon$\par\begin{equation}W^{\Dis\veps}(x,k)= \Dis\fa{1}{(2\pi)^{d}} \Dis\int_{\R^{d}} e^{iy.k} f^{\Dis\veps}(x-\Dis\veps y/2) \otimes f^{\Dis\veps \ast}(x+\Dis\veps y/2) dy\;.\end{equation}
\noindent\begin{proposition}\noindent Let the family  $f^{\Dis\veps}$ be uniformly bounded in $L^{2}(\R^{d})^n$. Then, upon extracting a subsequence, the semi 
classical transform $W^{\Dis\veps}$ converges weakly to a 
distribution $W(x,k)\in S^{\;'}(\R^{d}\times\R^{d})^n$, such 
that $Tr\ W(x,k)$ is a non-negative measure of bounded 
total mass (in the case $n=d$).\par
\end{proposition}
% Note that the positivity of the limit distribution together with the %identity %$$\Dis\int W_{\Dis\veps}(x,k) dk =|f_{\Dis\veps}(x)|^2 \;.$$\end{proposition}\noindent In the following, for simplicity, we consider the scalar case corresponding to $n=1$.\par
\noindent Let $a(x,k)$ be a test function in $S(\R^{d}\times\R^{d})$, where $x\in\R^{d}$ is the spatial variable, and $k\in\R^{d}$ is the momentum, or also the dual variable to $x$ in Fourier space. Then 
\begin{equation}<a,W^{\Dis\veps}>=(a^{w}(x,\Dis\varepsilon D) f^{\Dis\veps},f_{\Dis\veps})\end{equation}
\noindent where $<,>$ is the duality product between $S^{'}(\R^{d})$ and $S(\R^{d})$, $(,)$ is the $L^2(\R^{d})$ inner product, and the Weyl operator $a^{w}(x,\Dis\varepsilon D)$ is defined by 
\begin{equation}\left\{\begin{array}{ccc}[a^{w}(x,\Dis\varepsilon D)]f(x)= \Dis\fa{1}{(2\pi)^{d}}\Dis\int_{\R^{d}\times \R^{d}} a(\Dis\fa{x+y}{2},\Dis\varepsilon k) f(y) e^{i(x-y).k} dk dy \\[0.5cm]  
=\Dis\fa{1}{(2\pi\Dis\varepsilon )^{n}}\Dis\int_{\R^{d}} \hat{a}(\Dis\fa{x+y}{2},\Dis\fa{y-x}{\Dis\varepsilon}) f(y)dy\;.  \end{array}\right.\end{equation}
\noindent Here $\hat{a}$ is the Fourier transform of $a(x,k)$ in the variable $k $ only,  
\begin{equation}\hat{a}(x,y)=\Dis\int_{\R^{d}} e^{-ik.y} a(x,k) dk\;\end{equation}\noindent and this operator is bounded on $ L^{2}(\R^{d})$, uniformly in $\Dis\varepsilon $, 
\begin{equation}||a^{w}(x,\Dis\varepsilon D)||_{L^{2}(\R^{d})\rightarrow L^{2}(\R^{d})}\leq c(a)\;.    \end{equation}
\noindent We also introduce the pseudo differential operator at the scale $\varepsilon$, $a(x,\Dis\varepsilon D)$ by 
\begin{equation}[a(x,\Dis\varepsilon D)f](x)=\Dis\fa{1}{(2\pi)^{d}}\Dis\int_{\R^{d}} e^{ix.k}a(x,\Dis\varepsilon k) \hat{f}(k) dk\;.    \end{equation}
\noindent Again, one can show that the operators $a(x,\Dis\varepsilon D)$ are uniformly bounded on $L^{2}(\R^{d})$; there exists a constant $c(a)>0$ independent of $\Dis\varepsilon\in (0,1)$ (but depending on the function $a$) so that 
\begin{equation}||a(x,\Dis\varepsilon D)||_{L^{2}(\R^{d})\rightarrow L^{2}(\R^{d})}\leq c(a)\;\end{equation}
\noindent and furthermore, it satisfies for any $s>0$ \begin{equation}\label{esti1}\Dis\varepsilon^{s}||a(x,\Dis\varepsilon D)||_{H^{-s}(\R^{d})\rightarrow L^{2}(\R^{d})}\leq c_{s}(a)\;    \end{equation}
\noindent and
\begin{equation}\label{esti2}\Dis\varepsilon^{s}||a(x,\Dis\varepsilon D)||_{L^{2}(\R^{d})\rightarrow H^{s}(\R^{d})}\leq c_{s}(a)\;.    \end{equation}
\noindent The important point is that
\begin{equation}||a(x,\Dis\varepsilon D)-a^{w}(x,\Dis\varepsilon D)||_{L^{2}(\R^{d})\rightarrow L^{2}(\R^{d})}\longrightarrow 0\;\end{equation}as $\varepsilon \rightarrow 0$, so that the two quantizations are asymptotically equivalent.\par \noindent With the above notations, one has the following link between pseudo differential theory and semi classical transforms
\begin{equation}\Dis\lim_{\Dis\varepsilon \longrightarrow 0} (a(x,\Dis\varepsilon D)f^{\Dis\varepsilon},f^{\Dis\varepsilon})=<a,W>=Tr\Dis\int a(x,k) W(dx,dk)     \end{equation}
\noindent (where we have also included the vectorial case).\par\noindent We shall also need the following results, from pseudo differential calculus (adapted at the scale $\varepsilon$)
\noindent\begin{Lemma}\label{produ}\noindent The product of two operators  $a(x,\Dis\varepsilon D)$, $b(x,\Dis\varepsilon D)$ can be written as 
\begin{equation}\label{ppp}b(x,\Dis\varepsilon D) a(x,\Dis\varepsilon D)=(ba)(x,\Dis\varepsilon D)+\Dis\varepsilon/i (\nabla_{k} b.\nabla_{x} a)(x,\Dis\varepsilon D)+ \Dis\varepsilon^{2} Q_{\Dis\varepsilon}\end{equation}
\noindent where the operators $Q_{\Dis\varepsilon}$ are uniformly bounded on $L^{2}$ with respect to $\varepsilon$.\end{Lemma}
\noindent\begin{Lemma}{\bf (Localisation)} \label{locq} \noindent Let $f^{\Dis\varepsilon}(x)$ be a uniformly boundedfamiliy of functions in  $L^{2}$, and let $\mu_{f}(x,k)$ be any limit semi classical measure. Let $\phi(x)$ be a smooth function. Then the  semi classical  measure of the family $g^{\Dis\varepsilon}(x) =\phi(x) f^{\Dis\varepsilon}(x)$ is $|\phi(x)|^{2} \mu_{f}(x,k)$. Moreover, let $f^{\Dis\varepsilon}$, $g^{\Dis\varepsilon}$ be two uniformly bounded families of $L^{2}$ functions which coincide in an open neighbourhood of a point $x_{\;0}$. Then any limit semi classical  measure $\mu_{f}$  and  $\mu_{g}$ coincide in this neighbourhood.\par   \end{Lemma}
\smallskip
\section {\bf Proofs of the Theorems}\subsection{\bf Proof of Theorem (\ref{hsdaa}), Perfect boundary condition case}\noindent We consider the time harmonic form of Maxwell's system (\ref{3.4}) or equivalently (\ref{rajout1}), in the half space of $\R^{3}$, $(x^{3}\geq 0)$, where $E^{\Dis\varepsilon}=(E_{1}^{\Dis\varepsilon},E_{2}^{\Dis\varepsilon},E_{3}^{\Dis\varepsilon})$, and $x=(x',x^{3})$, $x'\in \R^{2}$, $x^3\geq 0$, with a perfect conductor boundary condition $\vec{n}\wedge \vec{E}^{\Dis\varepsilon}=0$, $n$ being the outward unit normal vector, i.e. $n=-\vec k$, which in our flat boundary case is equivalent to

\begin{equation}\left\{\begin{array}{ccccccc}E_{1}^{\Dis\varepsilon}=0\;,\\[0.4cm]E_{2}^{\Dis\varepsilon}=0\;.\end{array}\right.\end{equation}%\noindent {\bf Ca veut donc dire que les deux premieres composantes de ton vecteur $u^\varepsilon$ sont nulles sur le bord!! quelque chose dont tu ne tiens pas compte dans la suite me semble-t-il!!}\par
\noindent We set $E^{\Dis\varepsilon}$ to be zero in the lower half space $x^{3}<0$, and thus Maxwell's system (\ref{3.4}) or (\ref{rajout1}) can be rewritten as 
\begin{equation}\label{maxwel2}\fa{iw}{\Dis\varepsilon}A^{0}(x)u^{\displaystyle\varepsilon}+ \Dis\sum_{j=1}^{3} A^{j} \frac{\partial u^{\displaystyle\veps}}{\partial x_{j}}+Cu^{\displaystyle\veps}=f^{\displaystyle\veps}(x)+A_{b}u^{\displaystyle\veps}(x',0)\otimes\delta_{x_{3}=0}\;\end{equation}
\smallskip
\noindent where the "boundary" matrix $A_{b}$ is given by
\begin{equation}\label{bbbon}
A_{b}=\left(       \begin{array}{cccccccccccccccccccccc}  
       0 &  0  & 0 & 0 &0 &0 \\       0 &  0  & 0 & 0 &0 &0 \\       0 &  0  & 0 & 0 &0 &0 \\
       0 &  0  & 0 & 0 &-1 &0 \\
       0 &  0  & 0 & 1 &0 &0 \\
       0 &  0  & 0 & 0 &0 &0 \\
    \end{array}      \right)\;.
\end{equation}

\noindent In fact, let us recall that, \par
%\noindent {\bf Non, je ne suis pas d'accord jusqu' a la formule (3.40). Tout d'abord, parce que dans cette formule, tu melanges des vecteurs de taille 6 avec des vecteurs de taille 3!! donc c'est impossible. je t'ai appele au tel pour t'expliquer. Donc on est d'accord. Il faut partir des deux equations en E et H, prolonger et puis ecrire ca sous forme systeme hyperbolique. En conclusion, le membre de gauche de la formule ci-dessus ne change pas, mais par contre c'est le terme de bord qui va changer; la matrice $A^3$ qui apparait a droite change et on l'appelle maintenant $A_b$ pour au bord}.\par
%\noindent {\bf Au vu des calculs qu'on a fait par tel, je trouve que le terme de bord s'ecrit au signe pres comme suit: (j'ecris en ligne mais c'est des colonnes!!
%$$( 0 , 0 , 0 , -u_5 , u_4 , 0 )\delta_{x_3 =0}$$
%et donc le terme de bord devrait faire apparaitre la matrice $A_b$ avec $6$ zeros en ligne 1, 2 et 3. A la ligne 4 des zeros sauf en colonne 5 il y a -1. En ligne 5, des zeros sauf en colonne 4 ou il y a un 1, et enfin la ligne 6 avec des zeros!!}

\begin{equation}\left\{\begin{array}{ccccccc}\mbox{rot}E^{\Dis\varepsilon}=\nabla\wedge E^{\Dis\varepsilon}=\left( \begin{array}{cccccccccccccccccccccc}    \partial_{2} E_{3}^{\Dis\varepsilon}-\partial_{3} E_{2}^{\Dis\varepsilon}\\ \partial_{3} E_{1}^{\Dis\varepsilon}-\partial_{1} E_{3}^{\Dis\varepsilon}\\                    \partial_{1} E_{2}^{\Dis\varepsilon}-\partial_{3} E_{1}^{\Dis\varepsilon}\\ \end{array}    \right)=       \left(       \begin{array}{cccccccccccccccccccccc}          0 &  0  & 0 \\       0 &  0  & -1 \\       0 &  1  & 0  \\    \end{array}    \right)\partial_{1} E^{\Dis\varepsilon}+ \vspace{0.5cm}\\\hspace{1cm}\left(       \begin{array}{cccccccccccccccccccccc}          0 &  0  & 1 \\       0 &  0  & 0 \\      -1 &  0  & 0  \\    \end{array}    \right)\partial_{2}E^{\Dis\varepsilon}+\left(       \begin{array}{cccccccccccccccccccccc}          0 &  -1  & 0 \\       1 &  0  &  0 \\       0 &  0  &  0  \\    \end{array}    \right)\partial_{3} E^{\Dis\varepsilon}.\end{array}\right.\end{equation}
\noindent As $n\wedge E^{\Dis\varepsilon}=0$, we have that 
\begin{equation}\widetilde{\mbox{rot}E^{\Dis\varepsilon}}=\nabla\wedge \widetilde{E^{\Dis\varepsilon}}=\left( \begin{array}{cccccccccccccccccccccc}    \partial_{2} \tilde{E_{3}^{\Dis\varepsilon}}- \partial_{3} \tilde{E_{2}^{\Dis\varepsilon}}\\ \partial_{3} \tilde{E_{1}^{\Dis\varepsilon}}-\partial_{1} \tilde{E_{3}^{\Dis\varepsilon}}\\                    \partial_{1} \tilde{E_{2}^{\Dis\varepsilon}}-\partial_{3} \tilde{E_{1}^{\Dis\varepsilon}}\\ \end{array}   \right)\;.\end{equation}
\noindent For the magnetic field $H^{\Dis\varepsilon}$, let $\Omega \subseteq \R^{3}$ be a open domain de 
$\R^{3}$. Then $\forall \varphi\in C_{c}^{\infty}(\Omega)^3$, one has

\begin{equation}\Dis\int_{\Omega} \mbox{rot}H^{\Dis\varepsilon}.\varphi \;dx =
\Dis\int_{\Omega} H^{\Dis\varepsilon}.\mbox{rot}\varphi \;dx- 
\Dis\int_{\partial\Omega} (n\wedge H^{\Dis\varepsilon}).\varphi\; dx\;.
\end{equation}

\noindent Extending $ H^{\Dis\varepsilon}$ by zero in the full space, we have
\begin{equation}\label{abgna}\Dis\int_{\R^{3}} \widetilde{\mbox{rot}H^{\Dis\varepsilon}}.\varphi \;dx =
\Dis\int_{\R^{3}} \tilde{H}^{\Dis\varepsilon}.\mbox{rot}\varphi \;dx- 
\Dis\int_{\R^{3}} (n\wedge \tilde{H}^{\Dis\varepsilon}).\varphi\; dx\;.
\end{equation}

\noindent Here $\vec{n}=-\vec{k}$, thus we get
\begin{equation}\left\{\begin{array}{ccccccc}\widetilde{\mbox{rot}H^{\Dis\varepsilon}}=\left( \begin{array}{cccccccccccccccccccccc}    \partial_{2} \tilde{H_{3}^{\Dis\varepsilon}}- \partial_{3}\tilde{H_{2}^{\Dis\varepsilon}}\\ \partial_{3} \tilde {H_{1}^{\Dis\varepsilon}}-\partial_{1} \tilde{H_{3}^{\Dis\varepsilon}}\\                    \partial_{1} \tilde{H_{2}^{\Dis\varepsilon}}-\partial_{3} \tilde{H_{1}^{\Dis\varepsilon}}\\ \end{array}   \right)-\left(\begin{array}{cccccccccccccccccccccc}   H_{2}^{\Dis\varepsilon}\otimes \delta_{x_{3}=0}\\H_{1}^{\Dis\varepsilon}\otimes \delta_{x_{3}=0}\\0\\\end{array}   \right)=\vspace{0.5cm}\\\mbox{rot}\tilde{H^{\Dis\varepsilon}}- \left(\begin{array}{cccccccccccccccccccccc}   H_{2}^{\Dis\varepsilon}\otimes \delta_{x_{3}=0}\\H_{1}^{\Dis\varepsilon}\otimes \delta_{x_{3}=0}\\0\\\end{array} \right) = rot \tilde{H^{\Dis\varepsilon}} + \left(\begin{array}{cccccccccccccccccccccc}   H_{2}^{\Dis\varepsilon}\\H_{1}^{\Dis\varepsilon}\\0\\\end{array} \right) \otimes \delta_{x_{3}=0}\;\end{array}\right.\end{equation}
\noindent and using all the above notations, we get (\ref{maxwel2}).\par
\noindent Let $\theta(x)$ be a test function with compact support that is equal to one on a compact set $K$. We multiply $u^{\Dis\varepsilon}$ by $\theta(x)$, and thus we can define the semi classical  measure $\ddot\mu $ on $K$ for the family $\theta u^{\;\Dis\varepsilon}$, that we assume uniformly bounded in $L^2$.\par
\smallskip
\noindent More precisely, set    

\begin{equation}u^{\Dis\varepsilon ,\theta}(x)=\theta(x)u^{\Dis\varepsilon}(x)  \end{equation}
\noindent and let $u^{\Dis\varepsilon ,\theta}(x^{'},0)$ its boundary value, which is meaningfull in some negative Sobolev space, see \cite{chazarain} for instance. We shall assume 
that $u^{\Dis\varepsilon ,\theta}$ are uniformly 
bounded in $L^{2}(\R^2)$ and  that $u^{\Dis\varepsilon ,\theta}(x^{'},0) \delta_{x_3 =0}$ are uniformly 
bounded in $H^{-1/2-\alpha}(\R^3)$ (see \cite{cess} or \cite{chazarain}).\par
\noindent We let (after having possibly extracted a suitable sub-sequence) $\ddot\mu$ and $\ddot \nu$ be the (matrix valued) semi classical measures of $u^{\Dis\veps ,\theta}$ and $u^{\Dis\veps ,\theta} (x',0)$ resp.\par
\noindent Now Maxwell system can be rewritten, with the cuttof function $\theta$, as 
\begin{equation}\label{cut}\left\{\begin{array}{ccccccc}i\omega A^{0}(x) u^{\Dis\varepsilon ,\theta}+\Dis\varepsilon \sum_{j=1}^{3} A^{j}\Dis\fa{\partial u^{\Dis\varepsilon ,\theta}}{\partial x_{j}}-\Dis\varepsilon \sum_{j=1}^{3} A^{j}\Dis\fa{\partial \theta }{\partial x_{j}}u^{\Dis\varepsilon}(x)+\Dis\varepsilon C(x) u^{\Dis\varepsilon ,\theta}\\[0.5cm]=\Dis\varepsilon f^{\Dis\varepsilon ,\theta}(x)+\Dis\varepsilon A_bu^{\Dis\varepsilon ,\theta}(x^{'},0)\otimes \delta_{x_{3}=0}\;.\end{array}  \right.\end{equation}
\smallskip
\noindent Let $a(x,k)$ be a matrix-valued test function with compact support in $K$, with respect to $x$. Applying the operator $a_{\Dis\varepsilon }=a(x,\Dis\varepsilon D)$ on both sides of (\ref{cut}), and taking the inner product with $u^{\Dis\varepsilon ,\theta}$, we get 

\begin{equation}\label{helen}\left\{\begin{array}{ccccccc}(a_{\Dis\varepsilon }[i\omega A^{0}(x) u^{\Dis\varepsilon ,\theta}+\Dis\varepsilon \sum_{j=1}^{3} A^{j}\Dis\fa{\partial u^{\Dis\varepsilon ,\theta}}{\partial x_{j}}-\Dis\varepsilon \sum_{j=1}^{3} A^{j}\Dis\fa{\partial \theta }{\partial x_{j}}u^{\Dis\varepsilon}(x)+\Dis\varepsilon C(x) u^{\Dis\varepsilon ,\theta}],u^{\Dis\varepsilon ,\theta})\\[0.5cm]=\Dis\varepsilon (a_{\Dis\varepsilon}[f^{\Dis\varepsilon ,\theta}(x)],u^{\Dis\varepsilon ,\theta})+\Dis\varepsilon(a_{\Dis\varepsilon}[ A_bu^{\Dis\varepsilon ,\theta}(x^{'},0)\otimes\delta_{x_{3}=0}], u^{\Dis\varepsilon ,\theta})\;.\end{array}  \right.\end{equation}\noindent This is well defined in view of the usual rules of pseudo differential calculus (at the scale $\Dis\veps$).\par
%\noindent Recall that the matrix $A_{b}$ is a constant matrix given %by (\ref{bbbon}). Note that the 
%family %$A_bu_{\Dis\varepsilon}^{\theta}(x^{'},0)\otimes\delta_{x_{3}=0}$ is %uniformly bounded in $H^{-1/2-\alpha }$ for any $\alpha >0$.\par
\noindent To evaluate the limit of the second term of the right hand side in (\ref{helen}), let us set 
$$v^{\Dis\varepsilon}:=A_bu^{\Dis\varepsilon ,\theta}(x^{'},0)\otimes\delta_{x_{3}=0}$$ 
\noindent and thus we get that $v^{\Dis\varepsilon}$ is uniformly bounded in $H^{-1/2-\alpha }$ for any $\alpha >0$, with $\alpha<1/2$.\par
\noindent Next, note that \par
\begin{equation}\left\{\begin{array}{ccccccc}|\Dis\varepsilon (a_{\Dis\varepsilon }[v^{\Dis\varepsilon}],u^{\Dis\varepsilon ,\theta})|\leq \Dis\varepsilon | a_{\Dis\varepsilon }[v^{\Dis\varepsilon}]|\;||u^{\Dis\varepsilon ,\theta}|| _{L^{2}(\R^{3})} \leq \vspace{0.5cm}\\\Dis\varepsilon || a_{\Dis\varepsilon }||_{H^{-s}\longrightarrow L^{2}} ||v^{\Dis\varepsilon}||_{H^{-s}} ||u^{\Dis\varepsilon ,\theta}||_{L^{2}(\R^{3})}\leq c\fa{c_{s}}{\Dis\varepsilon^{s}}||v^{\Dis\varepsilon}||_{H^{-s}} ||u^{\Dis\varepsilon ,\theta}||_{L^{2}(\R^{3})}\;.\end{array}  \right.\end{equation}
\noindent Thus if we choose $s=-1/2-\alpha$, we get 
$$ ||\Dis\varepsilon (a_{\Dis\varepsilon }[v^{\Dis\varepsilon}],u^{\Dis\varepsilon ,\theta})||_{L^{2}(\R^{3})}\leq \Dis\fa{c \Dis\varepsilon^{1/2}}{\Dis\varepsilon^{\alpha}}\;,\;\;\; \alpha >0\;.$$
\noindent It follows that choosing $\alpha <\fa{1}{2}$, one has 
\begin{equation}\Dis\lim_{\Dis\varepsilon \longrightarrow 0}  \Dis\varepsilon (a_{\Dis\varepsilon }[v^{\Dis\varepsilon}],u^{\Dis\varepsilon ,\theta})=0\;.\end{equation}
\noindent Let us also show that the other terms in (\ref{helen}) are bounded (uniformly with respect to 
$\varepsilon$). Indeed, for the third term of (\ref{helen}), one has 

\begin{equation}\left\{\begin{array}{ccccccc}|(a_{\Dis\varepsilon }[-\Dis\varepsilon \sum_{j=1}^{3} A^{j}\Dis\fa{\partial \theta }{\partial x_{j}}u^{\Dis\varepsilon}],u^{\Dis\varepsilon ,\theta})|\leq \Dis\varepsilon | a_{\Dis\varepsilon }[\sum_{j=1}^{3} A^{j}\Dis\fa{\partial \theta }{\partial x_{j}}u^{\Dis\varepsilon}]|\; ||u^{\Dis\varepsilon ,\theta}||_{L^{2}(\R^{3})} \leq \vspace{0.4cm}\\\Dis\varepsilon || a_{\Dis\varepsilon }||_{L^{2}(\R^{3})} ||\sum_{j=1}^{3} A^{j}\Dis\fa{\partial \theta }{\partial x_{j}}u^{\Dis\varepsilon}||_{L^{2}(\R^{3})}||u^{\Dis\varepsilon ,\theta}||_{L^{2}(\R^{3})} \leq \Dis\varepsilon c\; 
\end{array}  \right.\end{equation}
\noindent and thus
\begin{equation}\Dis\lim_{\Dis\varepsilon \longrightarrow 0}  (-\Dis\varepsilon a_{\Dis\varepsilon }[\sum_{j=1}^{3} A^{j}\Dis\fa{\partial\theta}{\partial x_{j}}u^{\Dis\varepsilon}],u^{\Dis\varepsilon ,\theta})=0 \end{equation}\vspace{0.3cm}
\noindent and similary for the terms $\Dis\varepsilon(a_{\Dis\varepsilon }[C(x)u^{\Dis\varepsilon ,\theta}],u^{\Dis\varepsilon ,\theta})$, and $\Dis\varepsilon(a_{\Dis\varepsilon }[f^{\Dis\varepsilon ,\theta}],u^{\Dis\varepsilon ,\theta})$.\par
\noindent For the second term on the left hand side of (\ref{helen}), we set 
\begin{equation}\left\{\begin{array}{ccccccc}a_2:=\Dis\varepsilon (a_{\Dis\varepsilon }[ \sum_{j=1}^{3} A^{j}\Dis\fa{\partial u^{\Dis\varepsilon ,\theta}}{\partial x_{j}}],u^{\Dis\varepsilon ,\theta}) \;, \\[0.2cm]a_{\Dis\varepsilon }=a(x,\Dis\varepsilon D) \;, \\[0.2cm]b_{\Dis\varepsilon }= \Dis\varepsilon \sum_{j=1}^{3} A^{j}\Dis\fa{\partial } {\partial x_{j}} \end{array}  \right.\end{equation}
\noindent and using the product rule (\ref{ppp}), we get that 

\begin{equation}\left\{\begin{array}{ccccccc}a_2=(a_{\Dis\varepsilon }[b_{\Dis\varepsilon } u^{\Dis\varepsilon ,\theta},u^{\Dis\varepsilon ,\theta}) =(a(x,\Dis\varepsilon D)(b(x,\Dis\varepsilon D)[u^{\Dis\varepsilon ,\theta}],u^{\Dis\varepsilon ,\theta})\approx((a_{\Dis\varepsilon }b_{\Dis\varepsilon })[u^{\Dis\varepsilon ,\theta}],u^{\Dis\varepsilon ,\theta})\vspace{0.3cm}\\ \longrightarrow   \mbox{Tr} \Dis\int_{\R^3_{x}}\int_{\R^3_{k}} a(x,k) b(x,k) \ddot\mu(dx,dk)=- \mbox{Tr} \Dis\int_{\R^3_{x}}\int_{\R^3_{k}} a(x,k) (i\sum_{j=1}^{3} k_{j} A^{j})  \ddot\mu (dx,dk).\par \end{array}  \right.\end{equation}
\smallskip
\noindent For the first term of (\ref{helen}), we set 
\begin{equation}\left\{\begin{array}{ccccccc}a_1:=(a_{\Dis\varepsilon }[ i\omega A^{0}(x) u^{\Dis\varepsilon ,\theta}],u^{\Dis\varepsilon ,\theta}) \;, \\[0.2cm]a_{\Dis\varepsilon }=a(x,\Dis\varepsilon D)  \;, \\[0.2cm]b_{\Dis\varepsilon }=A^{0}(x) \end{array}  \right.\end{equation}
\noindent and thus, one has
\begin{equation}\left\{\begin{array}{ccccccc}a_1:=(a_{\Dis\varepsilon }[ i\omega b(x,\Dis\varepsilon D) u^{\Dis\varepsilon ,\theta}],u^{\Dis\varepsilon ,\theta})\approx i\omega(a_{\Dis\varepsilon }b_{\Dis\varepsilon }[u^{\Dis\varepsilon ,\theta}],u^{\Dis\varepsilon ,\theta})\vspace{0.2cm}\\ \longrightarrow i\omega \mbox{Tr} \Dis\int_{\R^3_{x}}\int_{\R^3_{k}} a(x,k) b(x,k) \ddot\mu (dx,dk)=i\omega \mbox{Tr} \Dis\int_{\R^3_{x}}\int_{\R^3_{k}} a(x,k) A^{0}(x) \ddot\mu (dx,dk) .\par \end{array}  \right.\end{equation}
\smallskip
\noindent Thus all in all, passing to the limit in (\ref{helen}), we get
\begin{equation}\mbox{Tr} \Dis\int_{\R^3_{x}}\int_{\R^3_{k}} a(x,k)[i\omega A^{0}-i\sum_{j=1}^{3} k_{j} A^{j}] \ddot\mu (dx,dk)=0\;\end{equation}
\noindent for all matrix valued test function $a(x,k)\in {\cal S}(\R^{3}\times\R^{3})$, which is equivalent to (in ${\cal S}'$)
\begin{equation}(i\omega A^{0}-i\sum_{j=1}^{3} k_{j} A^{j}) \ddot\mu =0\;\end{equation}
\noindent or
\begin{equation}\label{re}(-A^{0})[i \sum_{\Dis j=1}^{3}(A^{0})^{-1}k_{j}A^{j}-i\omega {\bf Id}] \ddot\mu =0\;\end{equation}
\noindent Let us set
\begin{equation}\label{gavuo}
L(x,k)=\sum_{\Dis j=1}^{3}(A^{0})^{-1}k_{j}A^{j}\;.
\end{equation}

\smallskip
\noindent In order to find the eigenvectors of the $6\times 6$ matrix $L(x,k )$, we shall use an orthonormal propagation basis of $\R^3$. We denote by $(\widehat{k},z^{1}(k),z^{2}(k))$ the orthonormal propagation triple consisting of the direction of propagation $\hat{k}=k/\mid k\mid$ and two transverse unit vectors $ z^{1}(k),z^{2}(k)$. In polar coordinates, they are, see for more details \cite{Geoleon}, \cite{Geoleonkel}
\begin{equation}\label{bas}\widehat{k}= {k\over {\mid k\mid}} =\left(\begin{array}{cccccccccccccccccccccc}   \sin \theta \cos \phi \\\sin \theta \sin \phi \\\cos \theta \end{array}  \right)\;,\hspace{0.2cm}z^{1}(k)=\left(\begin{array}{cccccccccccccc}   \cos\theta \cos \phi \\\cos \theta \sin \phi \\-\sin \theta \\\end{array}  \right)\;,\hspace{0.2cm}z^{2}(k)=\left(\begin{array}{ccccccccccccc}   -\sin \phi \\\cos \phi \\0\\ \end{array}  \right)\;\end{equation}
\noindent where $|k|=(k_{1}^{2}+k_{2}^{2}+k_{3}^{2})^{1/2}$.\par
\noindent Then the eigenvectors (which belong to $\R^6$) of the matrix $L(x,k )$ are given by [see \cite{Geoleon}, \cite{Geoleonkel}]
\begin{equation}\label{bbdde}
\left\{
\begin{array}{ccccccc}
b^1_0=\Dis\fa{1}{\sqrt \epsilon }(\hat{k},0) \;,\hspace{0.3cm} b^2_0=\Dis\fa{1}{\sqrt{\mu}}(0,\hat{k}) \vspace{0.2cm}\\
b_{+}^{1}=(\Dis\fa{1}{\sqrt{2\epsilon}} z^{1},\Dis\fa{1}{\sqrt{2\mu }}z^{2})\;,\hspace{0.3cm}
b_{+}^{2}=(\Dis\fa{1}{\sqrt{2\epsilon}} z^{2},-\Dis\fa{1}{\sqrt{2\mu}}z^{1})\vspace{0.2cm} \\   
b_{-}^{1}=(\Dis\fa{1}{\sqrt{2\epsilon}} z^{1},-\Dis\fa{1}{\sqrt{2\mu }}z^{2})\;,\hspace{0.3cm}
b_{-}^{2}=(\Dis\fa{1}{\sqrt{2\epsilon}} z^{2},\Dis\fa{1}{\sqrt{2\mu}}z^{1}).\\   
\end{array}  
\right.
\end{equation}
\noindent The eigenvectors $b^1_0$ and $b^2_0$ represent the non-propagating longitudinal and the other eigenvectors 
correspond to transverse modes of propagation with respect the speed of propagation $v$. These eigenvectors 
correspond to the eigenvalues listed in the following Lemma, whose proof follows from \cite{Geoleonkel}, \cite{Geoleon}.\smallskip
\noindent\begin{Lemma} 
\noindent The semi classical  measure $ \ddot{\mu}$ is supported on the set (recall that we assume that the frequency $\omega \neq 0$)\begin{equation}\label{ree}U= \left\{(x,k)\;,\;\omega_{+}=\omega \right\}\cup \left\{(x,k)\;,\; \omega_{-}=\omega \right\}\end{equation}\noindent where $v(x)=\Dis\fa{1}{\sqrt{\epsilon(x)\eta(x)}}\;$ is the propagation speed, and 
\smallskip
$\omega_{0}= \omega_0 (x,k) =0$ $\;,\; \omega_{+}= \omega_+ (x,k)=v(x)|k|$$\;,\; \omega_{-}=\omega_- (x,k)=-v(x)|k|$ 
are the eigenvalues (of constant multiplicity two) of the dispersion matrix $L$.\par
\end{Lemma}

\noindent Il follows that the semi classical measure $\ddot{\mu}(x,k)$ has the form 
\begin{equation}\label{decomposition}
\left\{
\begin{array}{ccccccc}\ddot{\mu}(x,k)= \mu_{+}^{1}(x,k) b_{+}^{1}(x,k)\otimes b_{+}^{1*}(x,k)+
\mu_{+}^{2}(x,k) b_{+}^{2}(x,k)\otimes b_{+}^{2*}(x,k)\vspace{0.2cm}\\
+\mu_{-}^{1}(x,k) b_{-}^{1}(x,k)\otimes b_{-}^{1*}(x,k)+
\mu_{-}^{2}(x,k) b_{-}^{2}(x,k)\otimes b_{-}^{2*}(x,k)
\end{array}  
\right.
\end{equation}
\noindent where $\mu_{+}^{1}\;,\mu_{+}^{2}$ are two scalar positive measures 
supported on the set $\left\{(x,k)\;,\;\omega_{+}=\omega \right\}$, and 
$\mu_{-}^{1}\;, \mu_{-}^{2}$, are two scalar positive measures 
supported on the set $\left\{(x,k)\;,\; \omega_{-}=\omega \right\}$. $b_{+}^{1}\;, 
b_{+}^{2}\;$ (resp. $ b_{-}^{1}\; b_{-}^{2}$) are the two eigenvectors of the matrix $L(x,k)$ given 
by (\ref{bbdde}), corresponding to the eigenvalue $\omega_{+}\;$ (resp. $ \omega_{-}$).\par
\noindent In view of the above reduction, we are led to find the transport equations for each of these four scalar semi classical  measures.\par
\noindent For this purpose, using the equation (\ref{cut}), we have the following identity 
\begin{center}\label{3.30}\begin{equation}\left\{\begin{array}{ccccccc}0=i\omega (a_{\Dis\varepsilon }[u^{\Dis\varepsilon ,\theta}],u^{\Dis\varepsilon ,\theta})-i\omega (a_{\Dis\varepsilon }[u^{\Dis\varepsilon ,\theta}],u^{\Dis\varepsilon ,\theta})=(i\omega a_{\Dis\varepsilon }[u^{\Dis\varepsilon ,\theta}],u^{\Dis\varepsilon ,\theta})+(a_{\Dis\varepsilon }[u^{\Dis\varepsilon ,\theta}],iwu^{\Dis\varepsilon ,\theta})\\[0.5cm]=(a_{\Dis\varepsilon }[-\Dis\varepsilon \sum_{\Dis j=1}^{3}(A^{0})^{-1}(x) A^{j} \Dis\fa{\partial u^{\Dis\varepsilon ,\theta}}{\partial x_{j}}+\Dis\varepsilon \sum_{\Dis j=1}^{3}(A^{0})^{-1}(x) A^{j} \Dis\fa{\partial \theta }{\partial x_{j}}u^{\Dis\varepsilon}(x) \\[0.5cm]-\Dis\varepsilon (A^{0})^{-1}(x)C(x) u^{\Dis\varepsilon ,\theta}+\Dis\varepsilon (A^{0})^{-1}(x)f^{\Dis\varepsilon ,\theta}(x)+\Dis\varepsilon (A^{0})^{-1}(x)  A_{b}u^{\Dis\varepsilon ,\theta}(x^{'},0)\otimes \delta_{x_{3}=0}],u^{\Dis\varepsilon ,\theta}) \vspace{0.3cm}\\+(a_{\Dis\varepsilon }[u^{\Dis\varepsilon ,\theta}],-\Dis\varepsilon \sum_{\Dis j=1}^{3}(A^{0})^{-1}(x) A^{j} \Dis\fa{\partial u^{\Dis\varepsilon ,\theta}}{\partial x_{j}}+\Dis\varepsilon \sum_{\Dis j=1}^{3}(A^{0})^{-1}(x) A^{j}\Dis\fa{\partial \theta }{\partial x_{j}}u^{\Dis\varepsilon}(x) \\[0.5cm]-\Dis\varepsilon (A^{0})^{-1}(x)C(x) u^{\Dis\varepsilon ,\theta}+\Dis\varepsilon (A^{0})^{-1}(x)f^{\Dis\varepsilon ,\theta}(x)+\Dis\varepsilon (A^{0})^{-1}(x)  A_{b}u^{\Dis\varepsilon ,\theta}(x^{'},0)\otimes \delta_{x_{3}=0}) \vspace{0.3cm}\\=-\Dis\varepsilon (a_{\Dis\varepsilon }[\sum_{\Dis j=1}^{3}(A^{0})^{-1}(x) A^{j} \Dis\fa{\partial u^{\Dis\varepsilon ,\theta}}{\partial x_{j}}],u^{\Dis\varepsilon ,\theta})+\Dis\varepsilon (\sum_{\Dis j=1}^{3} \Dis\fa{\partial }{\partial x_{j}}[(A^{0})^{-1}(x) A^{j}a_{\Dis\varepsilon } u^{\Dis\varepsilon ,\theta}],u^{\Dis\varepsilon ,\theta}) \\[0.5cm]+\Dis\varepsilon (a_{\Dis\varepsilon }[\sum_{\Dis j=1}^{3}(A^{0})^{-1}(x) A^{j}\Dis\fa{\partial \theta }{\partial x_{j}}u^{\Dis\varepsilon}(x)],u^{\Dis\varepsilon})- \Dis\varepsilon(a_{\Dis\varepsilon }[(A^{0})^{-1}(x) A_{b}u^{\Dis\varepsilon ,\theta}(x^{'},0) \otimes \delta_{x_{3}=0}],u^{\Dis\varepsilon ,\theta}) \\[0.5cm]+\Dis\varepsilon (a_{\Dis\varepsilon}[u^{\Dis\varepsilon ,\theta}],\sum_{\Dis j=1}^{3}(A^{0})^{-1}(x) A^{j} \Dis\fa{\partial \theta }{\partial x_{j}} u^{\Dis\varepsilon})-\Dis\varepsilon (a_{\Dis\varepsilon} [u^{\Dis\varepsilon ,\theta}],(A^{0})^{-1}(x)C(x) u^{\Dis\varepsilon ,\theta}],u^{\Dis\varepsilon ,\theta})\vspace{0.3cm}\\+\Dis\varepsilon (a_{\Dis\varepsilon }[u^{\Dis\varepsilon ,\theta}],(A^{0})^{-1}(x) A_{b}u^{\Dis\varepsilon ,\theta}(x^{'},0)\otimes \delta_{x_{3}=0})\;.\end{array}  \right.\end{equation}\end{center} 
\noindent Recalling that the function $\theta $ is equals to one identically on the support of $a(x,k)$, the third, fifth and sixth terms vanish at the limit, and thus the last equation can be rewritten as\par
\begin{equation}\label{lmi}\left\{\begin{array}{ccccccc}\Dis\varepsilon(a_{\Dis\varepsilon}[\sum_{\Dis j=1}^{3}(A^{0})^{-1}(x) A^{j}\Dis\fa{\partial u^{\Dis\varepsilon ,\theta}}{\partial x_{j}}],u^{\Dis\varepsilon ,\theta})+
\Dis\varepsilon (\sum_{\Dis j=1}^{3}\Dis\fa{\partial }{\partial x_{j}}[(A^{0})^{-1}(x) A^{j}a_{\Dis\varepsilon } u^{\Dis\varepsilon ,\theta}],u^{\Dis\varepsilon ,\theta}) \\[0.5cm]=\Dis\varepsilon (a_{\Dis\varepsilon}[(A^{0})^{-1}(x) A_{b}u^{\Dis\varepsilon ,\theta}(x^{'},0)\otimes \delta_{x_{3}=0}],u^{\Dis\varepsilon ,\theta})\\[0.5cm]
+\Dis\varepsilon (a_{\Dis\varepsilon }[u^{\Dis\varepsilon ,\theta}],(A^{0})^{-1}(x) A_{b}u^{\Dis\varepsilon ,\theta}(x^{'},0)\otimes \delta_{x_{3}=0})\;.\end{array}  \right.\end{equation}
\noindent Using the product rule (\ref{ppp}), it follows that  
\begin{equation}\label{prcc}\Dis\varepsilon \Dis\sum_{j=1}^{3} a_{\Dis\varepsilon } (A^{0})^{-1}(x) A^{j} \Dis\fa{\partial } {\partial x_{j}}-\Dis\varepsilon \Dis\sum_{j=1}^{3} \Dis\fa{\partial } {\partial x_{j}}A^{j} (A^{0})^{-1} a_{\Dis\varepsilon }=\phi_{0}(x,\Dis\varepsilon D)+ \Dis\varepsilon \phi_{1}(x,\Dis\varepsilon D)+ \Dis\varepsilon^{2} R_{\Dis\varepsilon}    
\end{equation}
\noindent where $\phi_{0}\;\;, \phi_{1}$, are given by 
\begin{equation}\left\{\begin{array}{ccccccc}\phi_{0}(x,k)=i a(x,k) \Dis\sum_{j=1}^{3} (A^{0})^{-1}(x)k_{j} A^{j}-i \Dis\sum_{j=1}^{3} k_{j} A^{j}(A^{0})^{-1}(x)a(x,k)\;,\\[0.5cm]\phi_{1}(x,k)= \Dis\sum_{j,m=1}^{3} \Dis\fa{\partial a} {\partial k_{m}}  \Dis\fa{\partial (A^{0})^{-1}}{\partial x^{m}}A^{j} k_{j}-\Dis\sum_{j=1}^{3} A^{j}\Dis\fa{\partial (A^{0})^{-1}}{\partial x^{j}} a- \Dis\sum_{j=1}^{3} A^{j} (A^{0})^{-1} \Dis\fa{\partial a} {\partial x^{j}}   \end{array}  \right.\end{equation}
\noindent and the operators $ R_{\Dis\varepsilon}$ are uniformly bounded on $L^{2}$.\par
\smallskip
\noindent On one hand, using (\ref{prcc}), and the two relations (\ref{esti1}), (\ref{esti2}), we pass 
to the limit in (\ref{lmi}), as $ \Dis\varepsilon \to 0$, and obtain
\begin{equation} \mbox{Tr}\int_{R^3_{x}}\int_{R^3_{k}} \phi_{0}(x,k) a(x,k)\ddot\mu(dx,dk)=0, \hspace{0.2cm} \forall a(x,k)     \end{equation}
\noindent which is already a known result (localization principle).\par
\noindent On the other hand, dividing (\ref{lmi}) by $\Dis\varepsilon $, and passing to the limit as $\Dis\varepsilon \to 0 $, we get 
\begin{equation}\label{lmimi} \mbox{Tr}\int_{\R^3_{x}}\int_{\R^3_{k}} \phi_{1}(x,k) \ddot\mu(dx,dk)+\Dis\lim_{\Dis\varepsilon\to 0} \Dis\fa{1}{\Dis\varepsilon} (\phi_{0}(x,\Dis\varepsilon D)u^{\Dis\varepsilon ,\theta},u^{\Dis\varepsilon ,\theta})= \Dis\lim_{\Dis\varepsilon\to 0} B_{\Dis\varepsilon}(a)\end{equation}
\noindent where 
\begin{equation}\label{termlim}B_{\Dis\varepsilon}(a)\equiv (a_{\Dis\varepsilon }[(A^{0})^{-1}(x) A_{b}u^{\Dis\varepsilon ,\theta} \otimes \delta_{x_{3}=0}],u^{\Dis\varepsilon ,\theta})+(a_{\Dis\varepsilon}[u^{\Dis\varepsilon ,\theta}],(A^{0})^{-1}(x) A_{b}u^{\Dis\varepsilon ,\theta}(x^{'})\otimes \delta_{x_{3}=0})\;.\end{equation}
\noindent In order to find the transport equation for the semi classical scalar measure $\ddot{\mu}$, we use the orthonormal propagation basis, and we consider first a test function $a(x,k)$ of the form
\begin{equation}\label{specia}a(x,k)=a_{+}(x,k) d_{+}^{1}(x,k)\otimes d_{+}^{1*}(x,k)\end{equation}
\noindent where $a_{+}(x,k)$ is any scalar smooth function, and 
\begin{equation}\label{tse}d_{+}^{1}(x,k)=A^{0}(x) b_{+}^{1}(x,k).  \end{equation}
%\noindent where $ b_{+}^{1}$ is the left eigenvector of the %dispersion matrix $L(x,k)$, corresponding to 
%the eigenvalue $\omega_{+}$.\par\noindent This choice is done to find the first positive scalar measure $\mu_+^1$ and we will then proceed similarly for the other measures. 
\noindent Recalling that $A^{0}$ is a symmetric matrix, with the choice (\ref{specia}), we note then 
that $\phi_{0}$ vanishes, while $\phi_{1}$ becomes 
\begin{equation}\label{cpcp}\left\{\begin{array}{ccccccc}\phi_{1}(x,k)= \Dis\sum_{m=1}^{3} \Dis\fa{\partial a_{+}(x,k)}{\partial k_{m}} d_{+}^{1}(x,k) \otimes d_{+}^{1*}(x,k)\Dis\fa{\partial (A^{0})^{-1}}{\partial x^{m}}\Dis\sum_{j=1}^{3} A^{j} k_{j}\\[0.3cm]-\Dis\sum_{j=1}^{3}\Dis\fa{\partial a_{+}(x,k)}{\partial x_{j}} A^{j}(A^{0})^{-1}
 d_{+}^{1}(x,k)\otimes d_{+}^{1*}(x,k) \\[0.3cm]+ a_{+}\{\Dis\sum_{m=1}^{3} \Dis\fa{\partial d_{+}^{1}} {\partial k_{m}}\otimes d_{+}^{1*} \Dis\fa{\partial (A^{0})^{-1}}
{\partial x^{m}}\Dis\sum_{j=1}^{3} A^{j} k_{j}+d_{+}^{1}(x,k)\\[0.3cm]
\otimes\Dis\sum_{m=1}^{3} \Dis\fa{\partial d_{+}^{1*}(x,k)}{\partial k_{m}} \Dis\fa{\partial (A^{0})^{-1}}{\partial x^{m}}\Dis\sum_{j=1}^{3} A^{j} k_{j}\\
-\Dis\sum_{j=1}^{3} A^{j}\Dis\fa{\partial (A^{0})^{-1}}{\partial x^{j}}d_{+}^{1}(x,k)\otimes  d_{+}^{1*}(x,k)\\[0.3cm]- \Dis\sum_{j=1}^{3} A^{j} (A^{0})^{-1}\Dis\fa{\partial d_{+}^{1}} {\partial x^{j}}\otimes d_{+}^{1*}(x,k)-\Dis\sum_{j=1}^{3} A^{j}(A^{0})^{-1} d_{+}^{1}(x,k)\otimes \Dis\fa{\partial d_{+}^{1*}} {\partial x^{j}}\}= 
\phi_{11}+ \phi_{12}+ \phi_{13}\;.  \end{array}  \right.\end{equation}
\smallskip 
\noindent where 
\begin{equation}\label{cpcpcp}\left\{\begin{array}{ccccccc}

\phi_{11}=\Dis\sum_{m=1}^{3} \Dis\fa{\partial a_{+}(x,k)}{\partial k_{m}} d_{+}^{1}(x,k) \otimes d_{+}^{1*}(x,k)\Dis\fa{\partial (A^{0})^{-1}}{\partial x^{m}}\Dis\sum_{j=1}^{3} A^{j} k_{j}\;, \\[0.3cm]\phi_{12}=-\Dis\sum_{j=1}^{3}\Dis\fa{\partial a_{+}(x,k)}{\partial x_{j}} A^{j}(A^{0})^{-1}
 d_{+}^{1}(x,k)\otimes d_{+}^{1*}(x,k) \;\;, \\[0.3cm]
\phi_{13}=a_{+}\{\Dis\sum_{m=1}^{3} \Dis\fa{\partial d_{+}^{1}} {\partial k_{m}}\otimes d_{+}^{1*} \Dis\fa{\partial (A^{0})^{-1}}
{\partial x^{m}}\Dis\sum_{j=1}^{3} A^{j} k_{j}+d_{+}^{1}(x,k)\\[0.3cm]
\otimes\Dis\sum_{m=1}^{3} \Dis\fa{\partial d_{+}^{1*}(x,k)}{\partial k_{m}} \Dis\fa{\partial (A^{0})^{-1}}{\partial x^{m}}\Dis\sum_{j=1}^{3} A^{j} k_{j}\\-\Dis\sum_{j=1}^{3} A^{j}\Dis\fa{\partial (A^{0})^{-1}}{\partial x^{j}}d_{+}^{1}(x,k)\otimes  d_{+}^{1*}(x,k)\\[0.3cm]- \Dis\sum_{j=1}^{3} A^{j} (A^{0})^{-1}\Dis\fa{\partial d_{+}^{1}} {\partial x^{j}}\otimes d_{+}^{1*}(x,k)-\Dis\sum_{j=1}^{3} A^{j}(A^{0})^{-1} d_{+}^{1}(x,k)\otimes \Dis\fa{\partial d_{+}^{1*}} {\partial x^{j}}\}\;.
\end{array}  \right.\end{equation}\noindent  We shall use the eigenvectors in the orthonormal basis (\ref{bas}), 
and the following normalization relations, 

\begin{equation}\label{renorm}\left\{\begin{array}{ccccccc}(A^{0}b_{\alpha},b_{\beta})=\delta_{\alpha \beta}\;,\\[0.5cm](A^{j} b_{+}, b_{+})=v \hat{k}_{j}.    \end{array}  \right.\end{equation}
\noindent We can then evaluate the first term  $\phi_{11}$ in (\ref{cpcp}). Indeed, we have 
\begin{equation}<\phi_{11},\ddot{\mu} >=< \Dis\fa{\partial a_{+}(x,k)}{\partial k_{m}} d_{+}^{1}(x,k)\otimes  d_{+}^{1*}(x,k)\Dis\sum_{j=1}^{3}\Dis\fa{\partial (A^{0})^{-1}}{\partial x^{j}}A^{j} k_{j},\mu_{+}^{1}(x,k) b_{+}^{1}(x,k) \otimes (b_{+}(x,k))^{1*})>\;.\end{equation}
\noindent Since 
\begin{equation}\Dis\fa{\partial }{\partial x^{m}}(\Dis\sum_{j=1}^{3} (A^{0})^{-1} k_{j}A^{j} d_{+}^{1})= \Dis\sum_{j=1}^{3} \Dis\fa{\partial (A^{0})^{-1}} {\partial x^{m}} k_{j}A^{j} d_{+}^{+1}+ \Dis\sum_{j=1}^{3} (A^{0})^{-1} k_{j}A^{j} \Dis\fa{\partial d_{+}^{1}}{\partial x^{m}}\;.  \end{equation}
\noindent It follows that 
\begin{equation}\Dis\sum_{j=1}^{3} \Dis\fa{\partial }{\partial x^{m}}((A^{0})^{-1} k_{j}A^{j} d_{+}^{1})-\Dis\sum_{j=1}^{3} (A^{0})^{-1} k_{j}A^{j} \Dis\fa{\partial d_{+}^{1}}{\partial x^{m}}= \Dis\sum_{j=1}^{3} \Dis\fa{\partial (A^{0})^{-1}} {\partial x^{m}} k_{j}A^{j} d_{+}^{1}\;. \end{equation}
\noindent Using the eigenvectors of the dispersion matrix in the orthonormal basis (\ref{bas}), one has 

\begin{equation}\Dis\sum_{j=1}^{3} \Dis\fa{\partial (A^{0})^{-1}} {\partial x^{m}} k_{j}A^{j}d_{+}^{1}=\Dis\fa{\partial \omega_{+}}{\partial x_{m}}d_{+}^{1}+ \omega_{+} \Dis\fa{\partial d_{+}^{1}}{\partial x_{m}}-\Dis\sum_{j=1}^{3} (A^{0})^{-1} k_{j}A^{j}\Dis\fa{\partial d_{+}^{1}}{\partial x^{m}}\;.\end{equation}
\noindent Thus the first term in (\ref{cpcp}), becomes 
\begin{equation}\left\{\begin{array}{ccccccc}<\phi_{11},\ddot{\mu}>=< \Dis\fa{\partial a_{+}}{\partial k_{m}} (b_{+}^{1},d_{+}^{1})(\Dis\sum_{j=1}^{3}\Dis\fa{\partial (A^{0})^{-1}}{\partial x^{m}} k_{j}A^{j} d_{+}^{1}  ,b_{+}^{1}),\mu_{+}^{1} > \\[0.3cm]=< \Dis\fa{\partial a_{+}}{\partial k_{m}} (\Dis\fa{\partial \omega_{+}}{\partial x^{m}}d_{+}^{1}+ \omega_{+} \Dis\fa{\partial d_{+}^{1}}{\partial x^{m}}- \Dis\sum_{j=1}^{3}(A^{0})^{-1} k_{j}A^{j}\Dis\fa{\partial d_{+}^{1}}{\partial x^{m}},b_{+}^{1}),\mu_{+}^{1}> \;.\end{array}  \right.\end{equation}
\noindent Using (\ref{renorm}), it follows that 
\begin{equation}\left\{\begin{array}{ccccccc}(\Dis\fa{\partial \omega_{+}}{\partial x^{m}}d_{+}^{1},b_{+}^{1})=(\Dis\fa{\partial \omega_{+}}{\partial x^{m}} A^{0}b_{+}^{1},b_{+}^{1})=\Dis\fa{\partial \omega_{+}}{\partial x^{m}} \;, m=1,2,3\\[0.4cm](\omega_{+} \Dis\fa{\partial d_{+}}{\partial x^{m}},b_{+}^{1})=\omega_{+}(\Dis\fa{\partial (A^{0}b_{+}^{1})}{\partial x^{m}},b_{+}^{1})=0\;, m=1,2,3\\[0.4cm](\Dis\sum_{j=1}^{3}(A^{0})^{-1} k_{j}A^{j}\Dis\fa{\partial d_{+}^{1}}{\partial x^{m}},b_{+}^{1}),b_{+}^{1})=(\Dis\sum_{j=1}^{3}(A^{0})^{-1} k_{j}A^{j}\Dis\fa{\partial ((A^{0})^{-1}b_{+}^{1})}{\partial x^{m}},b_{+}^{1})=0\;, m=1,2,3\;.  \end{array}  \right.\end{equation}
\noindent All in all, the first term of (\ref{cpcp}), becomes 

\begin{equation}\label{fist1}<\phi_{11},\ddot{\mu}>=<\Dis\sum_{m=1}^{3} \Dis\fa{\partial \omega_{+}}{\partial x^{m}}\Dis\fa{\partial a_{+}}{\partial k_{m}},\mu_{+}^{1} >\;.\end{equation}
\noindent For the second term $\phi_{12}$ in (\ref{cpcp}), one has 

\begin{equation}\left\{\begin{array}{ccccccc}<\phi_{12},\ddot{\mu}>=-< \Dis\sum_{j=1}^{3}\Dis\fa{\partial a_{+}(x,k)}{\partial x_{j}} (A^{0})^{-1} A^{j} d_{+}^{1}(x,k)\otimes d_{+}^{*1}(x,k),\mu_{+}^{1}> \\[0.4cm]=-<\Dis\sum_{j=1}^{3}\Dis\fa{\partial a_{+}(x,k)}{\partial x_{j}}(b_{+}^{1},d_{+}^{1}) ((A^{0})^{-1} A^{j}d_{+}^{1},b_{+}^{1}),\mu_{+}^{1} >.  \end{array}  \right. \end{equation}
\noindent Using the eigenvectors of the dispersion matrix in the orthonormal basis (\ref{bas}), 
and (\ref{renorm}), one has 
\begin{equation}\left\{\begin{array}{ccccccc}(b_{+}^{1},d_{+}^{1})=(b_{+}^{1},A^{0} b_{+}^{1})=1  \\[0.4cm]((A^{0})^{-1} A^{j}d_{+}^{1},b_{+}^{1})= ((A^{0})^{-1} A^{j} A^{0}b_{+}^{1},b_{+}^{1})=(A^{j}b_{+}^{1},b_{+}^{1})=\Dis\fa{\partial \omega_{+}}{\partial k_{j}}\end{array}  \right. \end{equation} 
\noindent and the second term $\phi_{12}$ in (\ref{cpcp}), becomes 

\begin{equation}\label{termnu2}<\phi_{12},\ddot{\mu} >=-< \Dis\sum_{j=1}^{3}\Dis\fa{\partial \omega_{+}}{\partial k_{j}}\Dis\fa{\partial a_{+}}{\partial x^{j}},\mu_{+}^{1} >\;.\end{equation}
\noindent For the third term $\phi_{13}$ in (\ref{cpcp}), we shall show that 
\begin{equation} \label{termnu3}<\phi_{13},\ddot{\mu}>=0\;. \end{equation}
\noindent We need to deal with the following term 
\begin{equation}\label{tremme}\left\{\begin{array}{ccccccc}<a_{+}\{ \Dis\sum_{m=1}^{3}\Dis\fa{\partial d_{+}^{1}} {\partial k_{m}} \otimes d_{+}^{1*} \Dis\fa{\partial (A^{0})^{-1}}{\partial x^{m}}\Dis\sum_{j=1}^{3}A^{j} k_{j}+
d_{+}^{1}(x,k)
\\[0.3cm]\otimes \Dis\sum_{m=1}^{3}\Dis\fa{\partial d_{+}^{1*}(x,k)}{\partial k_{m}} 
\Dis\fa{\partial (A^{0})^{-1}}{\partial x^{m}}\Dis\sum_{j=1}^{3}A^{j} k_{j}\\
- \Dis\sum_{j=1}^{3}A^{j}\Dis\fa{\partial (A^{0})^{-1}}{\partial x^{j}}d_{+}^{1}(x,k)
\otimes d_{+}^{1*}(x,k)\\[0.3cm]- \Dis\sum_{j=1}^{3}A^{j} (A^{0})^{-1}\Dis\fa{\partial d_{+}^{1}} {\partial x^{j}}\otimes d_{+}^{1*}(x,k)-\Dis\sum_{j=1}^{3}A^{j}(A^{0})^{-1} d_{+}^{1}(x,k)\otimes \Dis\fa{\partial d_{+}^{1*}} {\partial x_{j}}\},\mu_{+}^{1}(x,k) d_{+}^{1}(x,k)\otimes d_{+}^{1*}(x,k)> \;.\end{array}  \right.\end{equation}
\noindent Set
\begin{equation}\label{thir}\left\{\begin{array}{ccccccc}T=\Dis\sum_{m=1}^{3}(b_{+}^{1},\Dis\fa{\partial d_{+}^{1}}{\partial k_{m}})(b_{+}^{1}, \Dis\sum_{j=1}^{3}k_{j} A^{j} \Dis\fa{\partial (A^{0})^{-1}}{\partial x^{m}}d_{+}^{1})+\Dis\sum_{m=1}^{3}(b_{+}^{1},d_{+}^{1})(b_{+}^{1}, \Dis\sum_{j=1}^{3} k_{j}A^{j} \Dis\fa{\partial (A^{0})^{-1}}{\partial x^{m}}\Dis\fa{\partial d_{+}^{1}} {\partial k_{m}})\\[0.3cm]-\Dis\sum_{j=1}^{3}(b_{+}^{1},A^{j} \Dis\fa{\partial (A^{0})^{-1}}{\partial x^{j}}d_{+}^{1})(b_{+}^{1},d_{+}^{1})-\Dis\sum_{j=1}^{3}(b_{+}^{1},A^{j}(A^{0})^{-1}\Dis\fa{\partial d_{+}^{1}} {\partial x^{j}})(b_{+}^{1},d_{+}^{1}) \\[0.3cm]-\Dis\sum_{j=1}^{3}(b_{+}^{1},A^{j}(A^{0})^{-1} d_{+}^{1})(b_{+}^{1},\Dis\fa{\partial d_{+}^{1}}{\partial x_{j}}) \;.  \end{array}  \right.\end{equation}
\noindent We can rewrite (\ref{tremme}) as
\begin{equation}<\phi_{13},\ddot{\mu}>=< a_{+}[T],\mu_{+}^{1}>\;.\end{equation}
\noindent For the first term in (\ref{thir}), we use (\ref{bas}) and (\ref{renorm}), to get 

\begin{equation}\left\{\begin{array}{ccccccc}\Dis\sum_{m=1}^{3}(b_{+}^{1},\Dis\fa{\partial d_{+}^{1}}{\partial k_{m}})(b_{+}^{1},\Dis\sum_{j=1}^{3} k_{j} A^{j} \Dis\fa{\partial (A^{0})^{-1}}{\partial x^{m}}d_{+}^{1})=\Dis\sum_{m=1}^{3}(b_{+}^{1}, \Dis\fa{\partial \left[(A^{0})^{-1} b_{+}^{1}\right]}{\partial k_{m}})(b_{+}^{1}, \Dis\sum_{j=1}^{3}k_{j} A^{j}\Dis\fa{\partial (A^{0})^{-1}}{\partial x^{m}}d_{+}^{1}) \\[0.3cm]= \Dis\sum_{m=1}^{3}\overline{(\Dis\fa{\partial \left[(A^{0})^{-1} b_{+}^{1}\right]}{\partial k_{m}},b_{+}^{1})}(b_{+}^{1}, \Dis\sum_{j=1}^{3}k_{j} A^{j}\Dis\fa{\partial (A^{0})^{-1}}{\partial x^{m}}d_{+}^{1})\\[0.3cm]=\Dis\sum_{m=1}^{3}\overline{(\Dis\fa{\partial }{\partial k_{m}}\left[(A^{0})^{-1} b_{+}\right],b_{+}^{1})}(b_{+}^{1}, \Dis\sum_{j=1}^{3}k_{j} A^{j}\Dis\fa{\partial (A^{0})^{-1}}{\partial x^{m}}d_{+}^{1})=0\;\end{array}  \right.\end{equation}
\noindent and thus the first term in (\ref{thir}) vanishes. For the last term in (\ref{thir}), we use again 
(\ref{renorm}) to get %{\bf sommation sur les j? dans la formule qui suit?}
\begin{equation}\left\{\begin{array}{ccccccc}\Dis\sum_{j=1}^{3}(b_{+}^{1},A^{j}(A^{0})^{-1} d_{+}^{1})(b_{+}^{1},\Dis\fa{\partial d_{+}^{1}}{\partial x^{j}})=\Dis\sum_{j=1}^{3}(b_{+}^{1},A^{j}(A^{0})^{-1} A^{0}b_{+}^{1}) (b_{+}^{1},\Dis\fa{\partial d_{+}^{1}}{\partial x_{j}}) \\[0.4cm]=\Dis\sum_{j=1}^{3}(b_{+}^{1},A^{j}b_{+}^{1})(b_{+}^{1},\Dis\fa{\partial d_{+}^{1}}{\partial x_{j}})=\overline{\Dis\sum_{j=1}^{3}(A^{j}b_{+}^{1},b_{+}^{1})} (b_{+}^{1},\Dis\fa{\partial d_{+}^{1}}{\partial x_{j}})= v\hat{k}_{j} (b_{+}^{1},\Dis\fa{\partial d_{+}^{1}}{\partial x_{j}}) \\[0.4cm]=\Dis\sum_{j=1}^{3}\Dis\fa{\partial \omega_{+}}{\partial k_{j}}(b_{+}^{1},\Dis\fa{\partial d_{+}^{1}}{\partial x^{j}})\;.\end{array}  \right.\end{equation}
\noindent Thus (\ref{thir}) becomes %{\bf sommation sur les j? indice m? dans la formule qui suit?}
\begin{equation}\label{thir1}T= (b_{+}^{1},\Dis\sum_{j=1}^{3}\Dis\sum_{m=1}^{3}k_{j} A^{j} \Dis\fa{\partial (A^{0})^{-1}}{\partial x^{m}}\Dis\fa{\partial d_{+}^{1}}{\partial k_{m}}- \Dis\sum_{j=1}^{3}A^{j} \Dis\fa{\partial (A^{0})^{-1}}{\partial x^{j}} d_{+}^{1}-\Dis\sum_{j=1}^{3}A^{j} (A^{0})^{-1} \Dis\fa{\partial d_{+}^{1}}{\partial x_{j}}- \Dis\sum_{j=1}^{3}\Dis\fa{\partial \omega_{+}}{\partial k_{j}}\Dis\fa{\partial d_{+}^{1}}{\partial x^{j}})\;.\end{equation}
\noindent For the last term in (\ref{thir1}), we use (\ref{bas}) and (\ref{renorm}) to get 
%{\bf sommation sur les j? dans la formule qui suit?}
\begin{equation}\label{reforrr1}\left\{\begin{array}{ccccccc}(b_{+}^{1}, \Dis\sum_{j=1}^{3}\Dis\fa{\partial \omega_{+}}{\partial k_{j}} \Dis\fa{\partial d_{+}^{1}}{\partial x^{j}})=-\Dis\sum_{j=1}^{3}\Dis\fa{\partial \omega_{+}}{\partial k_{j}}(b_{+}^{1},\Dis\fa{\partial d_{+}^{1}}{\partial x^{j}})=-\Dis\sum_{j=1}^{3}\Dis\fa{\partial \omega_{+}}{\partial k_{j}}(b_{+}^{1}, \Dis\fa{\partial \left[A^{0}b_{+}^{1}\right]}{\partial x^{j}}) \\[0.4cm]= -\Dis\sum_{j=1}^{3}\Dis\fa{\partial \omega_{+}}{\partial k_{j}}\left\{(b_{+}^{1},\Dis\fa{\partial A^{0}}{\partial x^{j}}b_{+}^{1})+(b_{+}^{1},A^{0}\Dis\fa{\partial b_{+}^{1}}{\partial x^{j}})\right\}=\Dis\sum_{j=1}^{3}-\Dis\fa{\partial \omega_{+}}{\partial k_{j}}(b_{+}^{1},A^{0}\Dis\fa{\partial b_{+}^{1}}{\partial x^{j}})\;. \end{array}  \right.\end{equation}
\noindent For the second and third terms in (\ref{thir1}), we use (\ref{bas}) and (\ref{renorm}) to get 
%{\bf sommation sur les j? dans la formule qui suit?}
\begin{equation}\label{refor0}\left\{\begin{array}{ccccccc}-(b_{+}^{1},\Dis\sum_{j=1}^{3}\left\{A^{j} \Dis\fa{\partial (A^{0})^{-1}}{\partial x^{j}} d_{+}^{1}+A^{j} (A^{0})^{-1} \Dis\fa{\partial d_{+}^{1}} {\partial x_{j}}\right\})= (b_{+}^{1},\Dis\sum_{j=1}^{3}A^{j}\Dis\fa{\partial [(A^{0})^{-1}d_{+}^{1}]}{\partial x_{j}})\\[0.4cm]=-(A^{j}b_{+}^{1},\Dis\fa{\partial d_{+}^{1}} {\partial x_{j}}) \;.\end{array}  \right.\end{equation}
\noindent For the first term in (\ref{thir1}), we use (\ref{bas}) and (\ref{renorm}) to get 
%{\bf sommation sur les j? indice m? dans la formule qui suit?}

\begin{equation}\label{refor1}\left\{\begin{array}{ccccccc}(b_{+}^{1},\Dis\sum_{j=1}^{3}\Dis\sum_{m=1}^{3}k_{j} A^{j} \Dis\fa{\partial (A^{0})^{-1}}{\partial x^{m}}\Dis\fa{\partial d_{+}^{1}}{\partial k_{m}})=(b_{+}^{1},\Dis\sum_{j=1}^{3}k_{j} A^{j}\Dis\sum_{m=1}^{3} \Dis\fa{\partial (A^{0})^{-1}}{\partial x^{m}}\Dis\fa{\partial (A^{0}b_{+}^{1})}{\partial k_{m}}) \\[0.4cm]=(b_{+}^{1},\Dis\sum_{j=1}^{3}k_{j} A^{j} \Dis\sum_{m=1}^{3}\Dis\fa{\partial (A^{0})^{-1}} {\partial x^{m}} A^{0}\Dis\fa{\partial b_{+}^{1}}{\partial k_{m}})\;.\end{array}  \right.\end{equation}
\noindent But
\begin{equation}\Dis\fa{\partial }{\partial x^{m}} ( (A^{0})^{-1} A^{0})=\Dis\fa{\partial (A^{0})^{-1}}{\partial x^{m}} A^{0}+(A^{0})^{-1}\Dis\fa{\partial A^{0}}{\partial x^{m}}\;, m=1,2,3\;. \end{equation}
\noindent Thus (\ref{refor1}) becomes 
%{\bf sommation sur les j? indice m? dans la formule qui suit?}
\begin{equation}\label{refor2}\left\{\begin{array}{ccccccc}(b_{+}^{1}, \Dis\sum_{j=1}^{3}\Dis\sum_{m=1}^{3}A^{j}k_{j}\Dis\fa{\partial (A^{0})^{-1}} {\partial x^{m}} A^{0}\Dis\fa{\partial b_{+}^{1}}{\partial k_{m}})=-(b_{+}^{1},\Dis\sum_{j=1}^{3}\Dis\sum_{m=1}^{3}A^{j}k_{j}(A^{0})^{-1}\Dis\fa{\partial A^{0}}{\partial x^{m}}\Dis\fa{\partial b_{+}^{1}}{\partial k_{m}})\\[0.4cm]=-(b_{+}^{1},\Dis\sum_{j=1}^{3}A^{j}k_{j}(A^{0})^{-1})(b_{+}^{1},\Dis\sum_{m=1}^{3}\Dis\fa{\partial A^{0}}{\partial x^{m}}\Dis\fa{\partial b_{+}^{1}}{\partial k_{m}}) =-(b_{+}^{1}, \omega_{+}b_{+}^{1}) (b_{+}^{1},\Dis\sum_{m=1}^{3}\Dis\fa{\partial A^{0}}{\partial x^{m}}\Dis\fa{\partial b_{+}^{1}}{\partial k_{m}})\\[0.4cm]=-\omega_{+}\Dis\sum_{m=1}^{3}(b_{+}^{1},\Dis\fa{\partial A^{0}}{\partial x^{m}}\Dis\fa{\partial b_{+}^{1}}{\partial k_{m}})\;.\end{array}  \right.\end{equation}
\noindent Thus all in all, (\ref{thir1}) becomes 
%{\bf sommation sur les j? dans la formule qui suit?}
\begin{equation}\label{thir2}T=-\omega_{+}(b_{+}^{1},\Dis\sum_{j=1}^{3}\Dis\fa{\partial A^{0}}{\partial x^{j}}\Dis\fa{\partial b_{+}^{1}}{\partial k_{j}}) -(A^{j}b_{+}^{1},\Dis\fa{\partial b_{+}^{1}}{\partial x^{j}})+\Dis\fa{\partial \omega_{+}}{\partial k_{j}}(b_{+}, A^{0} \Dis\fa{\partial b_{+}^{1}} {\partial x^{j}})\;.\end{equation}
\noindent For the second and third terms in (\ref{thir2}), using the fact that $ b_{+}^{1}$ an eigenvector of the dispersion matrix, one has  %{\bf la formule qui suit n'est pas coherente!! sommation sur les j?}
\begin{equation}\label{thir3}\left\{\begin{array}{ccccccc}\Dis\fa{\partial }{\partial k_{j}}\{\Dis\sum_{j=1}^{3}(A^{0})^{-1} k_{j} A^{j}b_{+}^{1}\}=\Dis\fa{\partial }{\partial k_{j}}[\Dis\sum_{j=1}^{3}(A^{0})^{-1} k_{j} A^{j}] b_{+}^{1}+(\Dis\sum_{j=1}^{3}(A^{0})^{-1} k_{j} A^{j})\Dis\fa{\partial b_{+}^{1}} {\partial k_{j}} \\[0.4cm]= \Dis\fa{\partial }{\partial k_{j}} (\omega_{+}b_{+}^{1})= \Dis\fa{\partial \omega_{+}}{\partial k_{j}} b_{+}^{1}+\omega_{+}\Dis\fa{\partial b_{+}^{1}}{\partial k_{j}}\; \end{array}  \right.\end{equation}
\noindent which implies that
$$ A^{j}b_{+}^{1}=\Dis\fa{\partial \omega_{+}}{\partial k_{j}} A^{0}b_{+}^{1}+A^{0}\omega_{+}\Dis\fa{\partial b_{+}^{1}}{\partial k_{j}}-\Dis\sum_{j=1}^{3}k_{j} A^{j}\Dis\fa{\partial b_{+}^{1}} {\partial k_{j}}\;$$ 
\noindent and thus second and third terms in (\ref{thir2}) become 
%{\bf sommation sur les j? indice m? dans la formule qui suit?}
\begin{equation}-(A^{j}b_{+}^{1},\Dis\fa{\partial b_{+}^{1}}{\partial x^{j}})+\Dis\fa{\partial \omega_{+}}{\partial k_{j}}(b_{+}^{1}, A^{0} \Dis\fa{\partial b_{+}^{1}} {\partial x^{j}})=(b_{+}^{1},\Dis\sum_{j=1}^{3}k_{j} A^{j}\Dis\fa{\partial b_{+}^{1}} {\partial k_{j}}-A^{0}\omega_{+}\Dis\fa{\partial b_{+}^{1}}{\partial k_{j}})\;. \end{equation}
\noindent Thus all in all, we have 
%{\bf sommation sur les j? indice s? dans la formule qui suit?}
\begin{equation}\label{thir5}
\left\{\begin{array}{ccccccc}T=-\omega_{+}(b_{+}^{1},\Dis\sum_{j=1}^{3}\Dis\fa{\partial A^{0}}{\partial x^{j}}\Dis\fa{\partial b_{+}^{1}}{\partial k_{j}})+(\Dis\fa{\partial b_{+}^{1}}{\partial x^{m}},\Dis\sum_{j=1}^{3}k_{j} A^{j} \Dis\fa{\partial b_{+}^{1}}{\partial k_{m}}- A^{0}\omega_{+} \Dis\fa{\partial b_{+}^{1}}{\partial k_{m}})\\[0.4cm]
= \Dis\sum_{j=1}^{3}\Dis\fa{\partial \omega_{+}}{\partial x^{j}} (A^{0} b_{+}^{1}, \Dis\fa{\partial b_{+}^{1}}{\partial k_{j}})+0=0 \;\;\;, m=1,2,3
\end{array}  \right.\end{equation}
\noindent which yields (\ref{termnu3}).\par
%\noindent {\bf Enfin!! je n'ai bien sur pas verifie tes calculs!}\par
\smallskip
\noindent Now, using (\ref{fist1}), (\ref{termnu2}), (\ref{termnu3}), and integrating by parts, (\ref{lmimi}) becomes 

\begin{equation}\label{tarn}\left\{\begin{array}{ccccccc}<\phi_{11},\ddot{\mu}>+<\phi_{12},\ddot{\mu}>+<\phi_{13},\ddot{\mu}>=\Dis\sum_{m=1}^{3}< \Dis\fa{\partial \omega_{+}}{\partial x^{m}}\Dis\fa{\partial a_{+}}{\partial k_{m}},\mu_{+}^{1} >- \Dis\sum_{j=1}^{3}< \Dis\fa{\partial \omega_{+}}{\partial k_{j}}\Dis\fa{\partial a_{+}}{\partial x^{j}},\mu_{+}^{1} >\\[0.4cm]=<a_{+},\nabla_{k}\omega_{+} .\nabla_{x}\nu_{+}^{1} - \nabla_{x} \omega_{+}.\nabla_{k}\mu_{+}^{1} >=\Dis\lim_{ \Dis\varepsilon \to 0} B_{\Dis\varepsilon}(a)\;.\end{array}  \right.\end{equation}

\noindent There remains to determine the right hand side of (\ref{tarn}). Recall first that 

\begin{equation}\label {termlimll}B_{\Dis\varepsilon}(a)=(a_{\Dis\varepsilon }[(A^{0})^{-1}(x) A_{b}u^{\Dis\varepsilon ,\theta}(x',0) \otimes \delta_{x_{3}=0}],u^{\Dis\varepsilon ,\theta})+(a_{\Dis\varepsilon}[u^{\Dis\varepsilon ,\theta}],(A^{0})^{-1}(x) A_{b} u^{\Dis\varepsilon ,\theta}(x^{'},0)\otimes \delta_{x_{3}=0}).\end{equation}
\noindent Note that each term in (\ref{termlimll}) is of order $\Dis\varepsilon^{-1/2-\alpha}$ for any 
$\alpha > 0$ as can be seen from the $H^{s}$ estimates, since 
$u^{\Dis\varepsilon ,\theta}(x^{'},0)\otimes \delta_{x_{3}=0}$ is uniformly bounded in 
$H^{s}$ for $s=-1/2-\alpha $, for any $\alpha > 0$.\par
\noindent To get the limit of (\ref{termlimll}), we shall first use a special class of matrices $a(x,k)$ of the form 
\begin{equation}\label{clatest1}a(x,k)=\tilde{a}(x,k)[L(x,k)-\omega I]\;\end{equation}  
\noindent where $ L(x,k)$ is the dispersion matrix (\ref{gavuo}) and for any matrix $\tilde{a}(x,k)$ satisfying 
\begin{equation}\label{clatest2}\tilde{a}(x,k)[L(x,k)-\omega I]= [L^{*}(x,k)-\omega I]\tilde{a}(x,k)\;.\end{equation}
\noindent Using the test function (\ref{clatest1}) and the product rule (\ref{ppp}), the first term 
of (\ref{termlimll}) can be worked as follows 

\begin{equation}\label{fistrree}\left\{\begin{array}{ccccccc}\left(a_{\Dis\varepsilon }[ (A^{0})^{-1}(x) A_{b}u^{\Dis\varepsilon ,\theta} \otimes \delta_{x_{3}=0}],u^{\Dis\varepsilon ,\theta}\right)\\ [0.3cm] = \left(([L^{*}(x,k)-\omega I]\tilde{a}(x,\Dis\varepsilon D)\left[(A^{0})^{-1}(x) A_{b}u^{\Dis\varepsilon ,\theta} \otimes \delta_{x_{3}=0}\right],u^{\Dis\varepsilon ,\theta}\right) \\[0.3cm]\sim \left(([L^{*}(x,k)-\omega I](x,\Dis\varepsilon D)\tilde{a} (x,\Dis\varepsilon D)[(A^{0})^{-1}(x) A_{b}u^{\Dis\varepsilon ,\theta} \otimes \delta_{x_{3}=0}],u^{\Dis\varepsilon ,\theta}\right) \\[0.3cm] -\Dis\fa{\Dis\varepsilon}{i}\left((\nabla_{k}[L^{*}(x,k)-\omega I].\nabla_{x}\tilde{a})(x,\Dis\varepsilon D)[(A^{0})^{-1}(x) A_{b}u^{\Dis\varepsilon ,\theta} \otimes \delta_{x_{3}=0}],u^{\Dis\varepsilon ,\theta}\right)+\Dis\varepsilon^{2} \tilde Q_{\Dis\varepsilon}\\[0.3cm]\sim \left((\tilde{a} (x,\Dis\varepsilon D)[(A^{0})^{-1}(x) A_{b}u^{\Dis\varepsilon ,\theta} \otimes \delta_{x_{3}=0}],[L(x,k)-\omega I](x,\Dis\varepsilon D)u^{\Dis\varepsilon ,\theta}\right) \\[0.3cm] -\Dis\fa{\Dis\varepsilon}{i}\left((\nabla_{k}[L^{*}(x,k)-\omega I].\nabla_{x}\tilde{a})(x,\Dis\varepsilon D)[(A^{0})^{-1}(x) A_{b}u^{\Dis\varepsilon ,\theta} \otimes \delta_{x_{3}=0}],u^{\Dis\varepsilon ,\theta}\right)+\Dis\varepsilon^{2} \tilde Q_{\Dis\varepsilon}\;.\end{array}  \right.\end{equation}
with a term $\tilde Q_{\Dis\varepsilon}$ uniformly bounded.\par
\noindent The two last terms of the above formulae are uniformly bounded and vanishes to the limit.\par 
\noindent Indeed, for the first term, recall that 
\begin{equation}\label{eeeeg} i \Dis\sum_{j=1}^{3}A^{j}\Dis\fa{\partial u^{\Dis\varepsilon ,\theta}}{\partial x_{j}}=\omega A^{0}u^{\Dis\varepsilon ,\theta}\;.\end{equation}
\noindent We then use (\ref{helen}) and (\ref{eeeeg}), to rewrite (\ref{fistrree}) as 
\begin{equation}\label{fistrree1}\left\{\begin{array}{ccccccc}\left(a_{\Dis\varepsilon }[ (A^{0})^{-1}(x) A_{b}u^{\Dis\varepsilon ,\theta} \otimes \delta_{x_{3}=0}],u^{\Dis\varepsilon ,\theta}\right)\\[0.3cm] \sim (\tilde{a}(x,\Dis\varepsilon D)((A^{0})^{-1}(x) A_{b}u^{\Dis\varepsilon ,\theta}(x^{'})\otimes \delta_{x_{3}=0}),\Dis\fa{\Dis\varepsilon}{i}(A^{0})^{-1}(x)f^{\Dis\varepsilon ,\theta}(x)+\Dis\fa{\Dis\varepsilon}{i} (A^{0})^{-1}(x) A_{b}u^{\Dis\varepsilon ,\theta}(x^{'})\otimes \delta_{x_{3}=0} \\[0.3cm] +\Dis\fa{\Dis\varepsilon}{i} (A^{0})^{-1}(x) C(x)u^{\Dis\varepsilon ,\theta}(x)+\Dis\fa{\Dis\varepsilon}{i} \Dis\sum_{j=1}^{3}(A^{0})^{-1}(x) A^{j} \Dis\fa{\partial\theta }{\partial x^{j}}u^{\Dis\varepsilon}) \\[0.3cm]\sim (\tilde{a}(x,\Dis\varepsilon D)((A^{0})^{-1}(x) A_{b}u^{\Dis\varepsilon ,\theta}(x^{'})\otimes \delta_{x_{3}=0}),+\Dis\fa{\Dis\varepsilon}{i} (A^{0})^{-1}(x) A_{b}u^{\Dis\varepsilon ,\theta}(x^{'})\otimes \delta_{x_{3}=0} \\[0.3cm] +\Dis\fa{\Dis\varepsilon}{i} \Dis\sum_{j=1}^{3}(A^{0})^{-1}(x) A^{j} \Dis\fa{\partial\theta }{\partial x^{j}}u^{\Dis\varepsilon})\;. \\[0.3cm]\end{array}  \right.\end{equation}
\noindent For the second term of (\ref{termlimll}), still using the test function (\ref{clatest1}) 
and the product rule (\ref{ppp}), we get that 
\begin{equation}\label{fistrreeerr}\left\{\begin{array}{ccccccc}\left(a_{\Dis\varepsilon}[u^{\Dis\varepsilon ,\theta}],(A^{0})^{-1}(x) A_{b}u^{\Dis\varepsilon ,\theta}(x^{'})\otimes \delta_{x_{3}=0}\right)\\[0.3cm] =\left((\tilde{a}[L(x,k)-\omega I](x,\Dis\varepsilon D)[u^{\Dis\varepsilon ,\theta}],\left[(A^{0})^{-1}(x) A_{b}u^{\Dis\varepsilon ,\theta} \otimes \delta_{x_{3}=0}\right]\right) \\[0.3cm]\sim \left((\tilde{a}(x,\Dis\varepsilon D)\left[[L(x,k)-\omega I]u_{\Dis\varepsilon}^{\theta}\right],\left[(A^{0})^{-1}(x) A_{b}u^{\Dis\varepsilon ,\theta} \otimes \delta_{x_{3}=0}\right]\right) \\[0.3cm]-\Dis\fa{\Dis\varepsilon}{i}\left((\nabla_{k}\tilde{a}.\nabla_{x}[L(x,k)-\omega I])(x,\Dis\varepsilon D),\left[(A^{0})^{-1}(x) A_{b}u^{\Dis\varepsilon ,\theta} \otimes \delta_{x_{3}=0}\right]\right)+\Dis\varepsilon^{2} \tilde R_{\Dis\varepsilon}\;.\end{array}  \right.\end{equation}
with $\tilde R_{\Dis\varepsilon}$ uniformly bounded.\par
\noindent We use (\ref{helen}) and (\ref{eeeeg}) to rewrite (\ref{fistrreeerr}) as
\begin{equation}\label{fistrreeerrror}\left\{\begin{array}{ccccccc}\left(a_{\Dis\varepsilon}[u^{\Dis\varepsilon ,\theta}],(A^{0})^{-1}(x) A_{b}u^{\Dis\varepsilon ,\theta}(x^{'})\otimes \delta_{x_{3}=0}\right)=(\tilde{a}(x,\Dis\varepsilon D)( \Dis\fa{\Dis\varepsilon}{i} (A^{0})^{-1}(x) A_{b}u^{\Dis\varepsilon ,\theta}(x^{'})\otimes \delta_{x_{3}=0} \\[0.3cm] +\Dis\fa{\Dis\varepsilon}{i}(A^{0})^{-1}(x) A_{b}\Dis\fa{\partial\theta }{\partial x^{j}} u^{\Dis\varepsilon ,\theta} ,\left[(A^{0})^{-1}(x) A_{b}u^{\Dis\varepsilon ,\theta}(x^{'})\otimes \delta_{x_{3}=0}\right])\;.  \end{array}  \right.\end{equation}
\noindent Using these asymptotic expansions, passing to limit in (\ref{fistrree}), as $\Dis\varepsilon\to 0$, we obtain finally
\begin{equation}\lim_{ \Dis\varepsilon \to 0} B_{\Dis\varepsilon}(a)=0\;.\end{equation}

\noindent Now, we consider the general case of test functions in order to pass to the limit in the boundary term. For this purpose, we note that it is possible to write every test function $a_{+}$ as
\begin{equation}\label{classtest3}a_{+}(x,k)=a_{0}(x,k^{'})+a_{1}(x,k^{'})k_{3}+a_{2}(x,k)(v|k|-\omega)\end{equation}
\noindent where  $k=(k^{'},k_{3})$ and $ a_{0}\;,a_{1}$ and $a_{3}$ are scalar test functions, uniquely determined by $a_{+}$. For this point, we refer to \cite{GeoleonkelGuiller}. \par
\noindent In view of (\ref{classtest3}), we shall set

\begin{equation}\label{classtest3-1}
T_0 (a_+ ) =a_0 , \ T_1 (a_+ ) =a_1 \mbox{ and } T_2 (a_+ )=a_2.
\end{equation}
\noindent Then any $a$ of the form (\ref{specia}) can be written as 
\begin{equation}\label{testcla}\left\{\begin{array}{ccccccc}a(x,k)=(a_{0}(x,k^{'})+a_{1}(x,k^{'})k_{3})A^{0}(x) \\[0.3cm]+(a_{2}(x,k) d_{+}^{1}\otimes d_{+}^{1*}+ \Dis\fa{a_{0}(x,k^{'})+a_{1}(x,k^{'})k_{3}}{v|k|+\omega} d_{-}^{1}\otimes d_{-}^{1*}\\[0.3cm]+\Dis\fa{a_{0}(x,k^{'})+a_{1}(x,k^{'})k_{3}}{\omega} \Dis\sum_{j=1}^{2}d_{0}^{j}\otimes d_{0}^{j*})[L-\omega I](x,k)\;. \end{array}  \right. \end{equation}%\noindent {\bf Dans cette derniere formule, tu fais des sommes sur le j, mais il n'y a pas de j!!!}\par
\noindent Now, we note that the spectral representation of the matrix $L-\omega I $ can be written as 
%{\bf Ah, la formule qui vient est interessante!!}
\begin{equation}\label{rpopoa}L-\omega I=(\omega_{+}-\omega) b_{+}^{1}\otimes  d_{+}^{1*}+ (\omega_{-}-\omega) b_{-}^{1}\otimes  d_{-}^{1*}+(\omega_{0}-\omega) b^{0}\otimes d^{0*}\;.
\end{equation}
\noindent Recall that that the last term in (\ref{testcla}) has the same form of the test function $a=\tilde{a}[L-\omega I]$ of (\ref{clatest1}), and thus we can conclude for the limit of this term and we find that
\begin{equation}\lim_{ \Dis\varepsilon \to 0} B_{\Dis\varepsilon}(a_{2})=0\;.\end{equation}
\noindent Therefore, it is enough to find  the limits for $B_{\Dis\varepsilon}(a)$ only for the first two terms.\par
\noindent For this purpose, denote by  $a^{'}= a_{0}(x,k^{'})A^{0}(x)$ the first term in (\ref{testcla}). Multiplying it by a suitable cutoff function, $\phi(\Dis\varepsilon^{3} k_{3})$, with support compact, equal to one on a neighbourhood of zero, set
\begin{equation}a_{\Dis\varepsilon}^{''}=a^{'} \phi(\Dis\varepsilon^{3} k_{3})=[a_{0}(x,k^{'})A^{0}(x)]\phi(\Dis\varepsilon^{3} k_{3})\;.    \end{equation}
\noindent Using the product rule, the first term $B_{1\Dis\varepsilon}(a^{''})$ leads to 
\begin{equation}\label{testclaaa2}\left\{\begin{array}{ccccccc}B_{1\Dis\varepsilon}(a^{''})=B_{\Dis\varepsilon} \left(a_{\Dis\varepsilon}^{''} [ (A^{0})^{-1} A_{b}u_{\Dis\varepsilon}^{\theta}\otimes \delta_{x_3=0}],u^{\Dis\varepsilon ,\theta}\right) \\[0.3cm] \sim \left(a_{0}(x,\Dis\varepsilon D)A^{0}(x)\phi(\Dis\varepsilon^{3} k_{3})[ (A^{0})^{-1} A_{b}u^{\Dis\varepsilon ,\theta}\otimes \delta_{x_3=0}],u^{\Dis\varepsilon ,\theta}\right)\\[0.3cm]  \sim \Dis\int u^{\Dis\varepsilon ,\theta*}(x) dx \Dis\int \Dis\fa{d k}{(2\pi)^3} e^{ik.x}  a_{0}(x,\Dis\varepsilon k^{'}) \phi(\Dis\varepsilon^{3} k_{3})A_{b}\widehat{u^{\Dis\varepsilon ,\theta}(k^{'})}\end{array}  \right.\end{equation}
\noindent and similarly for the second term 
\begin{equation}\label{testclara}\left\{\begin{array}{ccccccc}B_{2\Dis\varepsilon}(a^{''})=B_{\Dis\varepsilon} \left(a_{\Dis\varepsilon}^{''} [u^{\Dis\varepsilon ,\theta}],(A^{0})^{-1} A_{b}u^{\Dis\varepsilon ,\theta}(x^{'})\otimes \delta_{x_{3}=0}]\right) \\[0.3cm] \Dis\int u^{\Dis\varepsilon ,\theta*}(x^{'}) dx^{'} \Dis\int \Dis\fa{d k}{(2\pi)^3} e^{ik^{'}.x^{'}}  a_{0}(x^{'},0,k^{'}) \phi(\Dis\varepsilon^{3} k_{3})A_{b}\widehat{u^{\Dis\varepsilon ,\theta}}(k) \;.\end{array}  \right.\end{equation}
\noindent Thus for the first term in (\ref{termlimll}), and for the test function written as (\ref{classtest3}), one has 

\begin{equation}\label{testclaooura}\Dis\lim_{\Dis\varepsilon \to 0}B_{\Dis\varepsilon}(a)=\Dis\lim_{\Dis\varepsilon \to 0}B_{1\Dis\varepsilon}(a_{0})+\Dis\lim_{\Dis\varepsilon \to 0}B_{2\Dis\varepsilon}(a_{0})=\mbox{Tr} \Dis\int A_{b} a_{0}(x^{'},0,k^{'}) d\ddot\nu \end{equation}
\noindent where $\ddot\nu$ is the semi classical measure of the boundary term  
$u_{\Dis\varepsilon}^{\theta}(x^{'},0)$.\par\noindent For the second term in (\ref{termlimll}), denoting by $a^{'}= a_{0}(x,k^{'}) k_{3}A^{0}(x)$, 
in the same way, we have 
\begin{equation}\label{testclaaa3}\left\{\begin{array}{ccccccc}B_{1\Dis\varepsilon}(a^{''})=B_{\Dis\varepsilon} \left(a_{\Dis\varepsilon}^{''} [ (A^{0})^{-1} A_{b}u^{\Dis\varepsilon ,\theta}\otimes \delta_{0}],u^{\Dis\varepsilon ,\theta}\right) \\[0.3cm] =\left(a_{1}(x,\Dis\varepsilon D)A^{0}(x)\phi(\Dis\varepsilon^{3} k_{3})[ (A^{0})^{-1} A_{b}u^{\Dis\varepsilon ,\theta}\otimes \delta_{x_3=0}],u^{\Dis\varepsilon ,\theta}\right)\\[0.3cm]  =\Dis\int u^{\Dis\varepsilon ,\theta*}(x) dx \Dis\int \Dis\fa{d k}{(2\pi)^3} e^{ik.x}  a_{1}(x,\Dis\varepsilon k^{'}) \Dis\varepsilon k_{3} \phi(\Dis\varepsilon^{3} k_{3})A_{b}\widehat{u^{\Dis\varepsilon ,\theta}(k^{'})}\;.\end{array}  \right.\end{equation}
\noindent Also 
\begin{equation}\label{testclaaarrra1}\left\{\begin{array}{ccccccc}B_{2\Dis\varepsilon}(a^{''})=B_{\Dis\varepsilon} \left(a_{\Dis\varepsilon}^{''} [u^{\Dis\varepsilon ,\theta}],(A^{0})^{-1} A_{b}u^{\Dis\varepsilon ,\theta}(x^{'})\otimes \delta_{x_{3}=0}]\right) \\[0.3cm] =\left(a_{1}(x,\Dis\varepsilon D) k_{3}A^{0}[u^{\Dis\varepsilon ,\theta}],(A^{0})^{-1} A_{b}u^{\Dis\varepsilon ,\theta}(x^{'})\otimes \delta_{x_{3}=0}]\right) \\[0.3cm] =\Dis\int u^{\Dis\varepsilon ,\theta*}(x^{'}) dx^{'} \Dis\int \Dis\fa{d k}{(2\pi)^3} e^{ik^{'}.x^{'}}  a_{1}(x^{'},0,k^{'}) \phi(\Dis\varepsilon^{3} k_{3})A_{b} \widehat{u^{\Dis\varepsilon ,\theta}}(k) \;.\end{array}  \right.\end{equation}
\noindent Thus, we have 
\begin{equation}\label{rentr}\left\{\begin{array}{ccccccc}B_{1\Dis\varepsilon}(a^{''})+B_{2\Dis\varepsilon}(a^{''})\sim -\Dis\fa{\Dis\varepsilon}{i} \Dis\int \Dis\fa{\partial {u}^{\Dis\varepsilon ,\theta*}}{\partial x^{3}}(x) dx \Dis\int \Dis\fa{d k}{(2\pi)^3} e^{ik.x}a_{1}(x,\Dis\varepsilon k^{'}) \phi(\Dis\varepsilon^{3} k_{3}) A_{b}\widehat{u^{\Dis\varepsilon ,\theta}}(k^{'}) \\[0.3cm]+\Dis\fa{\Dis\varepsilon}{i} \Dis\int u^{\Dis\varepsilon ,\theta*}(x) dx^{'} \Dis\int \Dis\fa{d k}{(2\pi)^3} e^{ik^{'}.x^{'}}  a_{1}(x^{'},0,k^{'})\phi(\Dis\varepsilon^{3} k_{3})A_{b}\Dis\fa{\widehat{\partial u^{\Dis\varepsilon ,\theta}}}{\partial x^{3}}(k)\;. \end{array}  \right.\end{equation}
\noindent Passing to limit in (\ref{rentr}), we get 

\begin{equation}\label{rez2}\Dis\lim_{\Dis\varepsilon \to 0}\left(B_{1\Dis\varepsilon}(a^{''})+B_{2\Dis\varepsilon}(a^{''})\right)=-\mbox{Tr} \Dis\int [\Dis\sum_{j=1}^{2} k_{j} A^{j}-\omega A^{0}(x^{'},0)]a_{1}(x^{'},0,k^{'}) d\ddot\nu  \;.  \end{equation}
\noindent Thus all in all, we get the limit of the boundary term (\ref{termlimll}) as 
\begin{equation}\label{testcla10}\Dis\lim_{\Dis\varepsilon \to 0} B_{\Dis\varepsilon}(a)=\mbox{Tr} \Dis\int [A_{b}(x^{'},0,k^{'})a_0 (x',0,k')-(\Dis\sum_{j=1}^{2} k_{j} A^{j}-\omega A^{0}(x^{'},0) )a_{1}(x^{'},0,k^{'}) ]d\ddot\nu \;.   \end{equation}
\noindent Note that if the test function $a(x,k)$ inside is supported away from $ x^{3}=0$ this limit equals zero, as it should be.\par
\noindent Now, because we are using special test functions satisfying (\ref{specia}), and since we are dealing only with $\mu^1_+$, we can as well assume that we are only seing the following part of $\ddot \nu$ given by\par$$\ddot \nu \sim \nu^1_{\alpha +} b^1_+ (k^+)\otimes b^{1\ast}_+ (k^+) +\nu^{1}_{\alpha\beta +} b^1_+ (k^+) \otimes b^{1\ast}_+ (k^-) +$$
$$+ \nu^{1}_{\beta\alpha +} b^1_+ (k^-) \otimes b^{1\ast}_+ (k^+) + \nu^{1}_{\beta +} b^1_+ (k^-) \otimes b^{1\ast}_+ (k^-).$$
\noindent This follows from the corresponding localization principle on the boundary.\par\noindent Next, note that 
$$\left[\Dis\sum_{j=1}^{2}(A^{0})^{-1} k_{j} A^{j}-\omega Id\right] b_{+}^{1}(k^{\pm})=
-k_{3}^{\pm} A^3b_{+}^{1}(k^{\pm})$$for scalar measures.\par
\noindent We have also
\begin{equation}\label{ercda}\left\{\begin{array}{ccccccc}(A_bb_{+}^{1}(k),b_{+}^{1}(k))= 0\;,\\[0.3cm]
(A^3b_{+}^{1}(k),b_{+}^{1}(k))=v\hat k^-_3 \;,\\[0.3cm](A^3b_{+}^{1}(k^{+}),b_{+}^{1}(k^{-}))=0\;,\\[0.3cm]
(A_bb_{+}^{1}(k^{+}),b_{+}^{1}(k^{-}))=0\;,\\[0.3cm](\Dis\sum_{j=1}^{2} k_{j} A^{j}-\omega A^{0})b_{+}^{1}(k^{\pm})=- k_{3}^{\pm} A^3b_{+}^{1}(k^{\pm})\;.\end{array}  \right.\end{equation}
\noindent Using (\ref{ercda}), (\ref{testcla}) and (\ref{classtest3}), the term (\ref{testcla10}) becomes
\begin{equation}\label{breba}\left\{\begin{array}{ccccccc}\Dis\lim_{ \Dis\varepsilon \to 0} B_{\Dis\varepsilon}(a)= \Dis\int \nu_{\alpha +}^{1} (dx^{'},dk^{'}) v {k}_{3}^{-} \widehat{k}_{3}^{-}  a_{1}  +\Dis\int \nu_{\beta +}^{1} (dx^{'},dk^{'}) v {k}_{3}^{+} \widehat{k}_{3}^{+} a_{1} \\[0.3cm]     \end{array}  \right.\end{equation}Recall that $v(x)=\Dis\fa{1}{\sqrt{\epsilon(x)\eta(x)}}$ is 
the propagation speed, the tangential vector $k'\in \R^{2}$, and the 
wave vector $k^{\pm}(k')=(k', k_{3}^{\pm})$ is defined by 
$$ k_{3}^{\pm}(x',0)=\pm\sqrt{\fa{\omega^{2}}{v(x',0)^{2}}-k'^{2}} \;.$$ 
\noindent By using formulas linked with the wave vectors, and in particuliar definitions given in (\ref{classtest3-1}), the above formulae reduces to

\begin{equation}\label{rajout2}
\left\{\begin{array}{ccccccc}\Dis\lim_{ \Dis\varepsilon \to 0} B_{\Dis\varepsilon}(a)=\\[0.3cm]
=\Dis\int \nu_{\alpha +}^{1} (dx^{'},dk^{'})  vk^-_3\hat k^-_3 T_1(a_+) (x',0,k')+\Dis\int \nu_{\beta +}^{1} (dx^{'},dk^{'}) 
vk^+_3\hat k^+_3 T_1(a_+) (x',0,k').
\end{array}
\right.
\end{equation}
%\\[0.3cm]
%=\Dis\int \nu_{\alpha}^{1} 
%(dx^{'},dk^{'})a_{+}(x^{'},0,k^{'},k_{3}^{-})v\hat{k}_{3}^{-} 
%+\Dis\int \nu_{\beta}^{1} %(dx^{'},dk^{'})a_{+}(x^{'},0,k^{'},k_{3}^{+})v\hat{k}_{3}^{+}\;.     \end{array}  \right.\end{equation}
\noindent  Combining (\ref{rajout2}) and (\ref{tarn}), we get the following  distributional form of the transport 
equation for the (scalar) positive measure $\nu_{+}^{1}(x,k)$  
\begin{equation}\label{transport}\nabla_{k}\omega_{+} .\nabla_{x}\mu^{1}_+ - \nabla_{x} \omega_{+}.\nabla_{k}\mu^{1}_+ = vk_3\hat k_3
[\nu^1_{\alpha +}T_1\delta_{k_3= k^-_3}  
+ \nu^1_{\beta +} T_1\delta_{k_3= k^+_3} ]\delta_{x_3 =0}.\end{equation}
\smallskip

\noindent  The other semi-classical measures in the formula are also dealt with in the same way as above, and we get 
\begin{equation}\label{transport-2}
\nabla_{k}\omega_{+} .\nabla_{x}\mu^{2}_+ - \nabla_{x} \omega_{+}.\nabla_{k}\mu^{2}_+ = vk_3\hat k_3
[\nu^2_{\alpha +}T_1 \delta_{k_3= k^-_3}  
+ \nu^2_{\beta +}T_1 \delta_{k_3=k^+_3} ]\delta_{x_3 =0}\\[0.3cm]
\end{equation}

\begin{equation}\label{transport-3}\nabla_{k}\omega_{+} .\nabla_{x}\mu^{1}_- - \nabla_{x} \omega_{+}.\nabla_{k}\mu^{1}_- = vk_3\hat k_3
[\nu^1_{\alpha -}T_1\delta_{k_3= k^-_3}  
+ \nu^1_{\beta -}T_1 \delta_{k_3= k^+_3} ]\delta_{x_3 =0}\end{equation}

\begin{equation}\label{transport-4}\nabla_{k}\omega_{+} .\nabla_{x}\mu^{2}_- - \nabla_{x} \omega_{+}.\nabla_{k}\mu^{2}_- = vk_3\hat k_3
[\nu^2_{\alpha -} T_1\delta_{k_3= k^-_3}  
+ \nu^2_{\beta -}T_1 \delta_{k_3=k^+_3} ]\delta_{x_3 =0}
\end{equation}

\subsection{\bf Proof of Theorem \ref{helen2}, Calderon type boundary condition}
\noindent In this case, for the exterior problem (\ref{exte}) (given in $\R^{3}_+$), extending by zero in the full space $\R^{3}$, we have the following eikonal equation for the exterior problem 
\begin{equation}\label{ccut45-exterieur}\left\{\begin{array}{ccccccc}i\omega A^{ext,0}(x) (u^{ext,\Dis\varepsilon ,\theta})+\Dis\varepsilon \sum_{j=1}^{3} A^{j}\Dis\fa{\partial (u^{ext,\Dis\varepsilon ,\theta})}{\partial x_{j}}-\Dis\varepsilon \sum_{j=1}^{3} A^{j}\Dis\fa{\partial \theta }{\partial x_{j}}u^{ext,\Dis\varepsilon}(x)
\\[0.5cm]+\Dis\varepsilon C^{ext}(x)(u^{ext,\Dis\varepsilon ,\theta}) =\Dis\varepsilon f^{\Dis\varepsilon ,\theta}(x)+\Dis\varepsilon A^3u^{ext,\Dis\varepsilon ,\theta}(x^{'},0)\otimes \delta_{x_{3}=0}.\end{array}  \right.\end{equation}
\noindent Note that, on the contrary of the perfect conductor case, we have not at this level taken into account Calderon transmission condition. We have also includede in the exterior field the incident one, using the same notations. Above the matrix $A^{ext,0}(x)$ is given by 
\begin{equation} A^{ext,0}  =\left(                  \begin{array}{cccccccccccccccccccccc}                      {\bf\epsilon}^{ext} {\bf Id}&  {\bf 0}  \\                    {\bf 0}      &  {\bf\eta}^{ext} {\bf Id} \\  \end{array}    \right)\end{equation}  
\noindent  where $ \epsilon^{ext} \;,\eta^{ext} $, are smoth functions in $C^{1}(\R^3) $, 
and the matrices $A^{\;j}$ are given by (\ref{elke}), and the matrix $ C^{ext} $ is given by
\begin{equation}\label{3.5}C^{ext}= \left(   \begin{array}{cccccccccccccccccccccc}       {\bf\sigma}^{ext}{\bf Id}& {\bf 0}\\     {\bf 0}     & {\bf 0} \\   \end{array}    \right)\end{equation}
\noindent with $\sigma^{ext}$ a smooth function in $C^{1}(\R^3) $, and
 $u^{int,\Dis\varepsilon ,\theta}(x^{'},0)$ is the boundary term for the 
interior problem (i.e. $x_{3}\leq 0$). In this case, we obtain that the dispersion matrix for 
the exterior problem is given by
\begin{equation}\label{matext}L^{ext}(x,k)=\Dis\sum_{j=1}^{3}(A^{ext,0})^{-1}k_{j}A^{j} \;.\end{equation}  
\noindent Recall that the matrix $L^{ext}$ has also three eigenvalues which constant multiplicity two. They are  
$$\omega_{0}^{ext}=0\hspace{0.2cm}, \omega_{+}^{ext}=v^{ext}|k^{'}|\hspace{0.2cm}, \omega_{-}^{ext}=-v^{ext}|k^{'}|$$ 
\noindent where $v^{ext}(x)=\Dis\fa{1}{\sqrt{\epsilon^{ext}(x)\eta^{ext}(x)}}$ is the propagation speed for the exterior problem.\par
\vspace{0.1cm}

\noindent As in the perfect conductor case, it follows that the associated semi classical measure $\ddot{\mu}^{ext}(x,k)$ has the form 
\begin{equation}\label{decomposition-exterieur}
\left\{
\begin{array}{ccccccc}\ddot{\mu}^{ext}(x,k)= \mu_{+}^{ext,1}(x,k) b_{+}^{ext,1}(x,k)\otimes b_{+}^{ext,1*}(x,k)+
\mu_{+}^{ext,2}(x,k) b_{+}^{ext,2}(x,k)\otimes b_{+}^{ext,2*}(x,k)\vspace{0.2cm}\\
+\mu_{-}^{ext,1}(x,k) b_{-}^{ext,1}(x,k)\otimes b_{-}^{ext,1*}(x,k)+
\mu_{-}^{ext,2}(x,k) b_{-}^{ext,2}(x,k)\otimes b_{-}^{ext,2*}(x,k)
\end{array}  
\right.
\end{equation}
\noindent where $\mu_{+}^{ext,1}\;,\mu_{+}^{ext,2}$ are two scalar positive measures 
supported on the set $\left\{(x,k)\;,\;\omega_{+}^{ext}=\omega \right\}$, and 
$\mu_{-}^{ext,1}\;, \mu_{-}^{ext,2}$, are two scalar positive measures 
supported on the set $\left\{(x,k)\;,\; \omega_{-}^{ext}=\omega \right\}$. $b_{+}^{ext,1}\;, 
b_{+}^{ext,2}\;$ (resp. $ b_{-}^{ext,1}\;,  b_{-}^{ext,2}$) are the two eigenvectors of the matrix $L^{ext}(x,k)$ given 
by (\ref{matext}), corresponding to the eigenvalue $\omega_{+}^{ext}\;$ (resp. $ \omega_{-}^{ext}$).\par
\noindent The semi classical  measure $ \ddot{\mu}^{ext}$ is supported on the set 
\begin{equation}\label{support-exterieur}U= \left\{(x,k)\;,\;\omega_{+}^{ext}=\omega \right\}\cup \left\{(x,k)\;,\; \omega_{-}^{ext}=\omega \right\}\;.\end{equation}

\noindent For instance, the transport equation for the first scalar measure is given by
\begin{equation}\label{transp-exterieur}\nabla_{k}\omega_{+}^{ext}.\nabla_{x}\mu^{ext,1}_+- \nabla_{x} \omega_{+}^{ext}.\nabla_{k}\mu^{ext,1}_+ =\\[0.3cm]v^{ext}\hat k_3[\nu^{ext,1}_{\alpha +} \delta_{k_3= k^{ext,-}_3}  
+ \nu^{ext,1}_{\beta +} \delta_{k_3= k^{ext,+}_3} ]\delta_{x_3 =0}\\[0.3cm]\end{equation}

\noindent where $\nu_{\alpha +}^{ext,1}\;, \nu_{\beta +}^{ext,1}$ are scalar measures corresponding to the boundary term $u^{ext,\Dis\varepsilon ,\theta}(x',0)$, and the 
wave vector $k^{ext,\pm}(k')=(k', k_{3}^{ext,\pm})$ is defined by 
$$ k_{3}^{ext,\pm}(x',0)=\pm\sqrt{\fa{\omega^{2}}{v^{ext}(x',0)^{2}}-k'^{2}} \;.$$ 
\noindent In fact, as in the previous sub-section, these (measures) coefficients come the following decomposition of the boundary semiclassical measure (seing only the first part of the set $U$)
$$\ddot \nu^{ext} \sim \nu^{ext,1}_{\alpha +} b^{ext,1}_+ (k^{ext,+})\otimes b^{ext,1\ast}_+ (k^{ext,+}) +\nu^{ext,1}_{\alpha\beta +} b^{ext,1}_+ (k^{ext,+}) \otimes b^{ext,1\ast}_+ (k^{ext,-}) +$$
$$+ \nu^{ext,1}_{\beta\alpha +} b^{ext,1}_+ (k^{ext,-}) \otimes b^{ext,1\ast}_+ (k^{ext,+}) + \nu^{ext,1}_{\beta +} b^{ext,1}_+ (k^{ext,-}) \otimes b^{ext,1\ast}_+ (k^{ext,-}).$$

\noindent For the interior problem (\ref{inte}), (given in $\R^{3}_-$), we have the following eikonal equation 
\begin{equation}\label{ccut45-interieur}\left\{\begin{array}{ccccccc}i\omega (A^{int,0})(x) (u^{int,\Dis\varepsilon ,\theta})+\Dis\varepsilon \sum_{j=1}^{3} A^{j}\Dis\fa{\partial (u^{int,\Dis\varepsilon ,\theta})}{\partial x_{j}}-\Dis\varepsilon \sum_{j=1}^{3} A^{j}\Dis\fa{\partial \theta }{\partial x_{j}}u^{int,\Dis\varepsilon}(x)
\\[0.5cm]+\Dis\varepsilon C^{int}(x)(u^{int,\Dis\varepsilon ,\theta})=\Dis\varepsilon f^{\Dis\varepsilon ,\theta}(x)+\Dis\varepsilon A^{int}_bu^{ext,\Dis\varepsilon ,\theta}(x^{'},0)\otimes \delta_{x_{3}=0}\;.\end{array}  \right.\end{equation}
\noindent where the matrix $A^{int,0}(x)$ is given by 
\begin{equation} A^{int,0}  =\left(                  \begin{array}{cccccccccccccccccccccc}                      {\bf \epsilon}^{int} {\bf Id}&  {\bf 0}  \\                     {\bf 0}      &  {\bf\eta}^{int} {\bf Id} \\  \end{array}    \right)\end{equation}  
\noindent where $ \epsilon^{int} \;,\eta^{int} $ are smooth functions in $C^{1}(\R^3) $, and 
the matrices $A^{\;j}$ are given by (\ref{elke}) and the matrix $ C^{int}$ is given by
\begin{equation}\label{3.5-1}C^{int}= \left(   \begin{array}{cccccccccccccccccccccc}       {\bf\sigma}^{int}{\bf Id}& {\bf 0}\\     {\bf 0}     & {\bf 0} \\   \end{array}    \right)\end{equation}\noindent with $\sigma^{int}$ a smooth function in $C^{1}(\R^3)$.\par

\noindent In this case, the interior dispersion matrix is given by
\begin{equation}\label{matint}L^{int}(x,k)=\Dis\sum_{j=1}^{3}(A^{int,0})^{-1}k_{j}A^{j} \end{equation}  
\noindent has three eigenvalues with constant multiplicity two. They are  
$$\omega_{0}^{int}=0\hspace{0.2cm}, \omega_{+}=v^{int}|k^{'}|\hspace{0.2cm}, \omega_{-}=-v^{int}|k^{'}|$$ 
\noindent where $v^{int}(x)=\Dis\fa{1}{\sqrt{\epsilon^{int}(x)\eta^{int}(x)}}$ is the propagation speed of the interior problem.\par
\vspace{0.1cm}

\noindent Again, it follows that the associated semi classical measure $\ddot{\mu}^{int}(x,k)$ has the form 
\begin{equation}\label{decomposition-interieur}
\left\{
\begin{array}{ccccccc}\ddot{\mu}^{int}(x,k)= \mu_{+}^{int,1}(x,k) b_{+}^{int,1}(x,k)\otimes b_{+}^{int,1*}(x,k)+
\mu_{+}^{int,2}(x,k) b_{+}^{int,2}(x,k)\otimes b_{+}^{int,2*}(x,k)\vspace{0.2cm}\\
+\mu_{-}^{int,1}(x,k) b_{-}^{int,1}(x,k)\otimes b_{-}^{int,1*}(x,k)+
\mu_{-}^{int,2}(x,k) b_{-}^{int,2}(x,k)\otimes b_{-}^{int,2*}(x,k)
\end{array}  
\right.
\end{equation}
\noindent where $\mu_{+}^{int,1}\;,\mu_{+}^{int,2}$ are two scalar positive measures 
supported on the set $\left\{(x,k)\;,\;\omega_{+}^{int}=\omega \right\}$, and 
$\mu_{-}^{int,1}\;, \mu_{-}^{int,2}$, are two scalar positive measures 
supported on the set $\left\{(x,k)\;,\; \omega_{-}^{int}=\omega \right\}$. $b_{+}^{int,1}\;, 
b_{+}^{int,2}\;$ (resp. $b_{-}^{int,1}\;,  b_{-}^{int,2}$) are the two eigenvectors of the matrix $L^{int}(x,k)$ given 
by (\ref{matint}), corresponding to the eigenvalue $\omega_{+}^{int}\;$ (resp. $ \omega_{-}^{int}$)s.\par
\noindent The semi classical  measure $ \ddot{\mu}^{int}$ is supported on the set 
\begin{equation}\label{support-interieur}U= \left\{(x,k)\;,\;\omega_{+}^{int}=\omega \right\}\cup \left\{(x,k)\;,\; \omega_{-}^{int}=\omega \right\}\;.\end{equation}

\noindent As an example, the transport equation for the first scalar measure is then
\begin{equation}\label{transport-interieur}\nabla_{k}\omega_{+}^{int}.\nabla_{x}\nu^{1int}_+- \nabla_{x} \omega_{+}^{int}.\nabla_{k}\nu^{1int}_+ =\\[0.3cm]v^{int}\hat{k}_{3}[\nu_{\alpha +}^{int,1}\delta_{k_{3}=k_{3}^{int,-}}+\nu_{\beta +}^{int,1}\delta_{k_{3}=k_{3}^{int,+}}]\delta_{x_3=0}\;.\end{equation}
\noindent where $\nu_{\alpha +}^{int,1}\;, \nu_{\beta +}^{int,1}$ are scalar measures associated with the semiclassical measure  
corresponding to the boundary term  $u^{int,\Dis\varepsilon ,\theta}(x',0)$, and the wave vector $k^{int,\pm}(k')=(k', k_{3}^{int,\pm})$ is defined by 
$$ k_{3}^{int,\pm}(x',0)=\pm\sqrt{\fa{\omega^{2}}{v^{int}(x',0)^{2}}-k'^{2}} \;.$$ 
\noindent Let us end by making some remarks about the scalar measures appearing the right hand sides of each transport equation, in the exterior as well as in the interior case.\par
\noindent Due to the Calderon transmission condition, it follows that one has the following, on the boundary $x_0 =0$
$$u^{int,\Dis\veps ,\theta} =M .u^{ext,\Dis\veps ,\theta}$$
where
$$M=\left(       \begin{array}{cccccccccccccccccccccc}  
      1 &  0  & 0 & 0 &0 &0 \\       0 &  1  & 0 & 0 &0 &0 \\       0 &  0  & 0 & 0 &0 &0 \\
       0 &  0  & 0 & 1 &0 &0 \\
       0 &  0  & 0 & 0 &1 &0 \\
       0 &  0  & 0 & 0 &0 &0 \\
    \end{array}      \right)\;.$$
\noindent It follows that ($M^2 =M$) 
$$\ddot \nu^{int} = M\ddot \nu^{ext}.$$
For instance, to get the scalar measure $\nu^{int,1}_{\alpha +}$, it is enough to take the trace of the above relation with 
$d^{int,1}_+ (k^{int,+})\otimes d^{int,1\ast}_+ (k^{int,+})$ (where we use left eigenvectors) and we get in this way
$$\nu^{int,1}_{\alpha +} =\ Tr (M\ddot \nu^{ext}.d^{int,1}_+ (k^{int,+})\otimes d^{int,1\ast}_+ (k^{int,+})).$$

\subsection{\bf Remarks on the curved interface case}\noindent In this case, we consider Maxwell's system above the surface given by $\Gamma:\ x_{3}=\phi(x^{'})$, where $x^{'}=(x_{1},x_{2})$, and $\phi\in W^{2}(\R^{2},\R)$ is a scalar function. 
\noindent We consider this system in time harmonic form, in the high frequency limit, and we consider a perfect boundary condition on $\Gamma $ . Again, we rewrite this system as a symetric one
\begin{equation}\label{maxwevbl25}\fa{iw}{\Dis\varepsilon}A^{0}(x)u_{\displaystyle\varepsilon}+ \Dis\sum_{j=0}^{3} A^{j} \frac{\partial u_{\displaystyle\veps}}{\partial x_{j}}+C u_{\displaystyle\veps}=f_{\displaystyle\veps}(x)+A^{3}u_{\displaystyle\veps}(x',0)\otimes \delta_{x_3 =\phi (x')}\end{equation}  
\noindent with $C(x)$ given by (\ref{3.522}), and $f_{\displaystyle\veps}\in L^{2}(\R^{3})^{3}$ and $A^{0}\;, A^{j}$ are given in (\ref{3.rr5}), (\ref{elke}).\par 
\noindent We shall reduce this curved case to a plane one, by introducing the new coordinates 
\begin{equation}y=x^{'}\;, z=x_{3}-\phi(x^{'})\;, \tilde x=(y,z). \end{equation}
\noindent Extending when necessary by zero in the all space $\R^{3}$, and thus (\ref{maxwevbl25}) becomes 
\begin{equation}\label{maxwel25}\fa{iw}{\Dis\varepsilon}\tilde A^{0}(y,z) v_{\displaystyle\varepsilon}(y,z)+ \Dis\sum_{j=0}^{3} \tilde A^{j} \frac{\partial v_{\displaystyle\veps}(y,z)}{\partial v}+\tilde C v_{\displaystyle\veps}(y,z)=\tilde f_{\displaystyle\veps}(y,z)+\tilde A_bv_{\displaystyle\veps}(y',0)\otimes\delta_{z=0}\;.\end{equation}  
\noindent where 
\begin{equation}\label{3.rrr5} \tilde A^{\;0} =\left(                  \begin{array}{cccccccccccccccccccccc}                      \epsilon(y,z) Id&  0  \\                     0      &  \eta(y,z) Id \\
  \end{array}    \right)\;,\tilde C= \left(   \begin{array}{cccccccccccccccccccccc}       \sigma(y,z)Id& 0\\    0     & 0 \\   \end{array}    \right)\end{equation}

\noindent are $3\times 3$ matrix valued smoth functions, with 
$ \epsilon(y,z) \;,\eta(y,z) \;, \sigma (y,z)$ smooth functions in $C^{1}(\R^3) $.\par

\noindent Set
\begin{equation}v_{\Dis\varepsilon}^{\theta}(y,z)=\theta(y,z)v_{\Dis\varepsilon}(y,z)  
\end{equation}
\noindent and the matrix of dispersion
\begin{equation}\label{matinter}L(\tilde x,k)=\Dis\sum_{j=1}^{3}((A^{\;0}))^{-1}k_{j}\tilde A^{j} \;\end{equation}   
\noindent where $\theta $ is a test function of compact support that is equal to one on a set compact $K$.\par 
\noindent Thus Maxwell system can be rewritten, with the cutoff function, as 
\begin{equation}\label{ccut}\left\{\begin{array}{ccccccc}i\omega \tilde A^{0}(y,z) v_{\Dis\varepsilon}^{\theta}(y,z)+\Dis\varepsilon \sum_{j=1}^{3}\tilde A^{j}\Dis\fa{\partial v_{\Dis\varepsilon}^{\theta}(y,z)}{\partial \tilde x_j}-\Dis\varepsilon \sum_{j=1}^{3}\tilde A^{j}\Dis\fa{\partial \theta (y,z)}{\partial \tilde x_j}v_{\Dis\varepsilon}(x)+\Dis\varepsilon \tilde C(y,z) v_{\Dis\varepsilon}^{\theta}(y,z)\\[0.5cm]=\Dis\varepsilon \tilde f_{\Dis\varepsilon}^{\theta}(y,z)+\Dis\varepsilon A_bv_{\Dis\varepsilon}^{\theta}(y,0)\otimes \delta_{z=0}\;.\end{array}  \right.\end{equation}
\noindent Then we can follow exactly the same steps as in the flat case.\par

\noindent { \bf Acknowledgements: }

\noindent The author would like to thank Radjesvarane Alexandre for several
discussions and suggestions during the preparation of this paper. \par


\begin{thebibliography}{0}


\bibitem{alexa} {\sc Alexandre.}
\newblock {\em Oscillations in PDE with singularities of Codimension One. Part I:
Review of the symbolic Calculus and Basic Definitions}. 
\newblock { Journal of Global Analysis and Geometry}. Preprint Submitted (2004).  


\bibitem{ant} {\sc Antonic, N.}
\newblock {\em H-measures applied to symmetric systems}.  
\newblock {\em Porc. Royal. Soc. Edinburgh 126A (1996), 1133-1155}. 



\bibitem{blp} {\sc Bensoussan, A., Lions, J.-L., Papanicolaou, G.}
\newblock {\em Asymptotic analysis for periodic structures}.  
\newblock North-Holland (1978).

\bibitem{belsi} {\sc Bleistein, N.}
\newblock {\em Mathematical Methods for wave phenomena,} 
\newblock {\em Academic Press,} (1984).    

\bibitem{bouix} {\sc Bouix, M.} {\em Les discontinuit\'es du rayonnement
\'electromagn\'etique}.  
\newblock Dunod, Paris (1965).


\bibitem{ngie} {\sc Burq N., Lebeau, G. } 
\newblock {\em Mesures de d\'efaut de compacit\'e, application au syst\`eme de Lam\'e}.
\newblock {\em Ann. Sci. Ecole Norm. Sup.}, (4)., 34, no. 6, 817-870, (2001).

\bibitem{chazarain} {\sc Chazarain, Piriou}
\newblock {\em Equations aux d\'eriv\'ees partielles lin\'eaires.}
\newblock {\em Presses du CNRS, Gauthier Villar, Paris (1985)}

\bibitem{cess} {\sc Cessenat, M.}
\newblock {\em Mathematical Methods In Electromagnetism. Linear Theory and
Applications}. 
\newblock World Scientific Publishing Co, Inc, River Edge, NJ (1996).

\bibitem{ciodon} {\sc Cioranescu, D., Donato, P.} 
\newblock {\em An introduction to homogenization}.  
\newblock Oxford University Press, New York Paris (1999).

   
\bibitem{Gér1} {\sc G\'erard, P.}    
\newblock {\em Microlocal defects measures, Comm.} PDEs, 16, (1991), 
1761-1794.


\bibitem{Gér2} {\sc G\'erard, P.} 
\newblock Mesures Semi-Classiques et Ondes de Bloch,
\newblock {\em S\'em. Ecole Polytechnique, expos\'e XVI,} (1990-91), 
1-19.


\bibitem{Gér3} {\sc G\'erard, P.} 
\newblock Oscillations and concentration effects in semilinear 
dispersive Wave equations,
\newblock {\em Journal of Functional Analysis,} 141, (1996), 60-98.

\bibitem{Géreric} {\sc G\'erard, P., Leichtnam., E.} 
\newblock {\em Ergodic properties of elegenfunctions for the Dirichlet problem}.
\newblock {\em Duke Mathematical Journal, 71 (2)} (1993), PP. 559-607. 

 
\bibitem{Gér2:Mar} {\sc G\'erard, P., Markovich, P., Mauser, N. and 
Poupaud, F.}  
\newblock {\em Homogenization limits and Wigner transforms},  
\newblock {\em Comm.Pure Appl. Math.,}  50, (1997), 323-380.

\bibitem{caclotid} {\sc Kammerer F.}  
\newblock {\em Propagation and absorption of concentration effects near shoch 
hypersurfaces for the heat equation}, 
\newblock {\em Asymptotic Analysis}, 
\newblock {\em 24., 107-141, (2000)}.

%\bibitem{hormr} {\sc Hormander, L.}
%\newblock {\em Analysis of linear partial differential operateurs III.} Pseudo-differential operators. Corrected reprint 
%of the 1985 original. Fundamental principles of Mathematical sciences. 274, Springer-verlay, Berlin, (1994).
  

\bibitem{lmul} {\sc Miller L.}
\newblock {\em Refraction of high-frequency waves density by sharp interfaces and semi-classical measures at 
the boundary}. Preprint., Ecole Polytechnique Palaiseau. (1999). 


\bibitem{ngde} {\sc N\'ed\'elec, J.C.} 
\newblock {\em Acoustic and electromagnetic equation. Integral representations for harmonic problems.,} 
Applied Mathematical Sciences, 144, Springer-verlag, New york, (2001).


\bibitem{Lipa} {P.L. Lions and T. Paul,}
\newblock Sur les Mesures de Wigner, Revista Mat. Iberoamericana, 9, 1993, 553-618. 

\bibitem{Geoleon} {\sc Papanicolaou, G., Ryzhik, L.} 
\newblock {\em Waves and transport}. In IAS/Park City Mathematics Series, Vol. 5., Caffarelli Weinan E, eds., AMS., (1998), 305-382.

\bibitem{Geoleonkel} {\sc Papanicolaou, G., Ryzhik, L., Keller, J.} 
\newblock {\em Transport equations for elastic and other waves in random media}. Waves Motion, 24., (1996), 327-370.

\bibitem{GeoleonkelGui} {\sc Papanicolaou, G., Keller., J.B., Bal., G., Ryzhik, L.} 
\newblock {\em Transport Theory for Acoustic waves with reflection and transmission at interfaces}. Wave Motion., 30, (1999), 303-327.

\bibitem{GeoleonkelGuilla} {\sc Papanicolaou, G., Keller., J.B., Bal., G., Ryzhik, L.} 
\newblock {\em Radiative Transport in a periodic structure}. Journal of Statistical physics., 95, (1/2)., 479-494, (1999).

\bibitem{GeoleonkelGuiller} {\sc Papanicolaou, G., Keller., J.B.,G., Ryzhik, L.} 
\newblock {\em Transport equations for waves in a half space}. Comm. PDE's, 22., (1997), 1869-1911.

\bibitem{poumark} {\sc Poupaud, F., Markowich, P.A.}
\newblock {\em The Maxwell equation in a periodic medium; homogenization of the 
energy density.} Ann. Scuola Norm. Sup. Pisa cl. Sci-(4), 23, no. 2, 301-324, (1996). 

\bibitem{tar1} {\sc Tartar, L.} 
\newblock {\em Cours Peccot}. Coll\`ege de France. Unpublished (1977).

\bibitem{tar2} {\sc Tartar, L.} 
\newblock {\em Compensated compactness and applications to partial differential equations}. Nonlinear analysis and mechanics, Vol.IV, 
Pitman, Boston, Mass-London., (1979), 136-212.
          
\bibitem{tar} {\sc Tartar, L.}
\newblock {\em H-measures, a new approach for studing homogenization, 
oscillations and concentration effects in partial differential equations.} Proc. Roy. Soc. Edinburgh., 115A, (1990), 193-230.

\bibitem{tay1} {\sc Taylor, M.E.}
\newblock {\em Partial differential equations: Basic theory.}
Springer-verlag, New york 563 pp (1996).

\end{thebibliography}
\end{document}